%
\newskip\ttglue
\font\fiverm=cmr5
\font\fivei=cmmi5
\font\fivesy=cmsy5
\font\fivebf=cmbx5
\font\sixrm=cmr6
\font\sixi=cmmi6
\font\sixsy=cmsy6
\font\sixbf=cmbx6
\font\sevenrm=cmr7
\font\eightrm=cmr8
\font\eighti=cmmi8
\font\eightsy=cmsy8
\font\eightit=cmti8
\font\eightsl=cmsl8
\font\eighttt=cmtt8
\font\eightbf=cmbx8
\font\ninerm=cmr9
\font\ninei=cmmi9
\font\ninesy=cmsy9
\font\nineit=cmti9
\font\ninesl=cmsl9
\font\ninett=cmtt9
\font\ninebf=cmbx9
\font\twelverm=cmr12
\font\twelvei=cmmi12
\font\twelvesy=cmsy12
\font\twelveit=cmti12
\font\twelvesl=cmsl12
\font\twelvett=cmtt12
\font\twelvebf=cmbx12


\def\eightpoint{\def\rm{\fam0\eightrm}  
  \textfont0=\eightrm \scriptfont0=\sixrm \scriptscriptfont0=\fiverm
  \textfont1=\eighti  \scriptfont1=\sixi  \scriptscriptfont1=\fivei
  \textfont2=\eightsy  \scriptfont2=\sixsy  \scriptscriptfont2=\fivesy
  \textfont3=\tenex  \scriptfont3=\tenex  \scriptscriptfont3=\tenex
  \textfont\itfam=\eightit  \def\it{\fam\itfam\eightit}
  \textfont\slfam=\eightsl  \def\sl{\fam\slfam\eightsl}
  \textfont\ttfam=\eighttt  \def\tt{\fam\ttfam\eighttt}
  \textfont\bffam=\eightbf  \scriptfont\bffam=\sixbf
    \scriptscriptfont\bffam=\fivebf  \def\bf{\fam\bffam\eightbf}
  \tt  \ttglue=.5em plus.25em minus.15em
  \normalbaselineskip=9pt
  \setbox\strutbox=\hbox{\vrule height7pt depth2pt width0pt}
  \let\sc=\sixrm  \let\big=\eightbig \normalbaselines\rm}

\def\eightbig#1{{\hbox{$\textfont0=\ninerm\textfont2=\ninesy
        \left#1\vbox to6.5pt{}\right.$}}}


\def\ninepoint{\def\rm{\fam0\ninerm}  
  \textfont0=\ninerm \scriptfont0=\sixrm \scriptscriptfont0=\fiverm
  \textfont1=\ninei  \scriptfont1=\sixi  \scriptscriptfont1=\fivei
  \textfont2=\ninesy  \scriptfont2=\sixsy  \scriptscriptfont2=\fivesy
  \textfont3=\tenex  \scriptfont3=\tenex  \scriptscriptfont3=\tenex
  \textfont\itfam=\nineit  \def\it{\fam\itfam\nineit}
  \textfont\slfam=\ninesl  \def\sl{\fam\slfam\ninesl}
  \textfont\ttfam=\ninett  \def\tt{\fam\ttfam\ninett}
  \textfont\bffam=\ninebf  \scriptfont\bffam=\sixbf
    \scriptscriptfont\bffam=\fivebf  \def\bf{\fam\bffam\ninebf}
  \tt  \ttglue=.5em plus.25em minus.15em
  \normalbaselineskip=11pt
  \setbox\strutbox=\hbox{\vrule height8pt depth3pt width0pt}
  \let\sc=\sevenrm  \let\big=\ninebig \normalbaselines\rm}

\def\ninebig#1{{\hbox{$\textfont0=\tenrm\textfont2=\tensy
        \left#1\vbox to7.25pt{}\right.$}}}


\def\twelvepoint{\def\rm{\fam0\twelverm}  
  \textfont0=\twelverm \scriptfont0=\eightrm \scriptscriptfont0=\sixrm
  \textfont1=\twelvei  \scriptfont1=\eighti  \scriptscriptfont1=\sixi
  \textfont2=\twelvesy  \scriptfont2=\eightsy  \scriptscriptfont2=\sixsy
  \textfont3=\tenex  \scriptfont3=\tenex  \scriptscriptfont3=\tenex
  \textfont\itfam=\twelveit  \def\it{\fam\itfam\twelveit}
  \textfont\slfam=\twelvesl  \def\sl{\fam\slfam\twelvesl}
  \textfont\ttfam=\twelvett  \def\tt{\fam\ttfam\twelvett}
  \textfont\bffam=\twelvebf  \scriptfont\bffam=\eightbf
    \scriptscriptfont\bffam=\sixbf  \def\bf{\fam\bffam\twelvebf}
  \tt  \ttglue=.5em plus.25em minus.15em
  \normalbaselineskip=11pt
  \setbox\strutbox=\hbox{\vrule height8pt depth3pt width0pt}
  \let\sc=\sevenrm  \let\big=\twelvebig \normalbaselines\rm}

\def\twelvebig#1{{\hbox{$\textfont0=\tenrm\textfont2=\tensy
        \left#1\vbox to7.25pt{}\right.$}}}
\magnification=\magstep1
\def\firstpage{1}
\pageno=\firstpage
\font\fiverm=cmr5
\font\sevenrm=cmr7
\font\sevenbf=cmbx7
\font\eightrm=cmr8
\font\eightbf=cmbx8
\font\ninerm=cmr9
\font\ninebf=cmbx9
\font\tenbf=cmbx10
\font\magtenbf=cmbx10 scaled\magstep1

\font\magnineeufm=eufm9 scaled\magstep1

\catcode`\@=11
%

\def\undefine#1{\let#1\undefined}
\def\newsymbol#1#2#3#4#5{\let\next@\relax
 \ifnum#2=\@ne\let\next@\msafam@\else
 \ifnum#2=\tw@\let\next@\msbfam@\fi\fi
 \mathchardef#1="#3\next@#4#5}
\def\mathhexbox@#1#2#3{\relax
 \ifmmode\mathpalette{}{\m@th\mathchar"#1#2#3}%
 \else\leavevmode\hbox{$\m@th\mathchar"#1#2#3$}\fi}
\def\hexnumber@#1{\ifcase#1 0\or 1\or 2\or 3\or 4\or 5\or 6\or 7\or 8\or
 9\or A\or B\or C\or D\or E\or F\fi}

\font\tenmsa=msam10
\font\sevenmsa=msam7
\font\fivemsa=msam5
\newfam\msafam
\textfont\msafam=\tenmsa
\scriptfont\msafam=\sevenmsa
\scriptscriptfont\msafam=\fivemsa
\edef\msafam@{\hexnumber@\msafam}
\mathchardef\dabar@"0\msafam@39
\def\dashrightarrow{\mathrel{\dabar@\dabar@\mathchar"0\msafam@4B}}
\def\dashleftarrow{\mathrel{\mathchar"0\msafam@4C\dabar@\dabar@}}

\def\ulcorner{\delimiter"4\msafam@70\msafam@70 }
\def\urcorner{\delimiter"5\msafam@71\msafam@71 }
\def\llcorner{\delimiter"4\msafam@78\msafam@78 }
\def\lrcorner{\delimiter"5\msafam@79\msafam@79 }
\def\yen{{\mathhexbox@\msafam@55}}
\def\checkmark{{\mathhexbox@\msafam@58}}
\def\circledR{{\mathhexbox@\msafam@72}}
\def\maltese{{\mathhexbox@\msafam@7A}}

\font\tenmsb=msbm10
\font\sevenmsb=msbm7
\font\fivemsb=msbm5
\newfam\msbfam
\textfont\msbfam=\tenmsb
\scriptfont\msbfam=\sevenmsb
\scriptscriptfont\msbfam=\fivemsb
\edef\msbfam@{\hexnumber@\msbfam}
\def\Bbb#1{{\fam\msbfam\relax#1}}
\def\widehat#1{\setbox\z@\hbox{$\m@th#1$}%
 \ifdim\wd\z@>\tw@ em\mathaccent"0\msbfam@5B{#1}%
 \else\mathaccent"0362{#1}\fi}
\def\widetilde#1{\setbox\z@\hbox{$\m@th#1$}%
 \ifdim\wd\z@>\tw@ em\mathaccent"0\msbfam@5D{#1}%
 \else\mathaccent"0365{#1}\fi}
\font\teneufm=eufm10
\font\seveneufm=eufm7
\font\fiveeufm=eufm5
\newfam\eufmfam
\textfont\eufmfam=\teneufm
\scriptfont\eufmfam=\seveneufm
\scriptscriptfont\eufmfam=\fiveeufm

\catcode`\@=11
\newsymbol\boxdot 1200
\newsymbol\boxplus 1201
\newsymbol\boxtimes 1202
\newsymbol\square 1003
\newsymbol\blacksquare 1004
\newsymbol\centerdot 1205
\newsymbol\lozenge 1006
\newsymbol\blacklozenge 1007
\newsymbol\circlearrowright 1308
\newsymbol\circlearrowleft 1309
\undefine\rightleftharpoons
\newsymbol\rightleftharpoons 130A
\newsymbol\leftrightharpoons 130B
\newsymbol\boxminus 120C
\newsymbol\Vdash 130D
\newsymbol\Vvdash 130E
\newsymbol\vDash 130F
\newsymbol\twoheadrightarrow 1310
\newsymbol\twoheadleftarrow 1311
\newsymbol\leftleftarrows 1312
\newsymbol\rightrightarrows 1313
\newsymbol\upuparrows 1314
\newsymbol\downdownarrows 1315
\newsymbol\upharpoonright 1316
 
\newsymbol\downharpoonright 1317
\newsymbol\upharpoonleft 1318
\newsymbol\downharpoonleft 1319
\newsymbol\rightarrowtail 131A
\newsymbol\leftarrowtail 131B
\newsymbol\leftrightarrows 131C
\newsymbol\rightleftarrows 131D
\newsymbol\Lsh 131E
\newsymbol\Rsh 131F
\newsymbol\rightsquigarrow 1320
\newsymbol\leftrightsquigarrow 1321
\newsymbol\looparrowleft 1322
\newsymbol\looparrowright 1323
\newsymbol\circeq 1324
\newsymbol\succsim 1325
\newsymbol\gtrsim 1326
\newsymbol\gtrapprox 1327
\newsymbol\multimap 1328
\newsymbol\therefore 1329
\newsymbol\because 132A
\newsymbol\doteqdot 132B
 
\newsymbol\triangleq 132C
\newsymbol\precsim 132D
\newsymbol\lesssim 132E
\newsymbol\lessapprox 132F
\newsymbol\eqslantless 1330
\newsymbol\eqslantgtr 1331
\newsymbol\curlyeqprec 1332
\newsymbol\curlyeqsucc 1333
\newsymbol\preccurlyeq 1334
\newsymbol\leqq 1335
\newsymbol\leqslant 1336
\newsymbol\lessgtr 1337
\newsymbol\backprime 1038
\newsymbol\risingdotseq 133A
\newsymbol\fallingdotseq 133B
\newsymbol\succcurlyeq 133C
\newsymbol\geqq 133D
\newsymbol\geqslant 133E
\newsymbol\gtrless 133F
\newsymbol\sqsubset 1340
\newsymbol\sqsupset 1341
\newsymbol\vartriangleright 1342
\newsymbol\vartriangleleft 1343
\newsymbol\trianglerighteq 1344
\newsymbol\trianglelefteq 1345
\newsymbol\bigstar 1046
\newsymbol\between 1347
\newsymbol\blacktriangledown 1048
\newsymbol\blacktriangleright 1349
\newsymbol\blacktriangleleft 134A
\newsymbol\vartriangle 134D
\newsymbol\blacktriangle 104E
\newsymbol\triangledown 104F
\newsymbol\eqcirc 1350
\newsymbol\lesseqgtr 1351
\newsymbol\gtreqless 1352
\newsymbol\lesseqqgtr 1353
\newsymbol\gtreqqless 1354
\newsymbol\Rrightarrow 1356
\newsymbol\Lleftarrow 1357
\newsymbol\veebar 1259
\newsymbol\barwedge 125A
\newsymbol\doublebarwedge 125B
\undefine\angle
\newsymbol\angle 105C
\newsymbol\measuredangle 105D
\newsymbol\sphericalangle 105E
\newsymbol\varpropto 135F
\newsymbol\smallsmile 1360
\newsymbol\smallfrown 1361
\newsymbol\Subset 1362
\newsymbol\Supset 1363
\newsymbol\Cup 1264
 
\newsymbol\Cap 1265
 
\newsymbol\curlywedge 1266
\newsymbol\curlyvee 1267
\newsymbol\leftthreetimes 1268
\newsymbol\rightthreetimes 1269
\newsymbol\subseteqq 136A
\newsymbol\supseteqq 136B
\newsymbol\bumpeq 136C
\newsymbol\Bumpeq 136D
\newsymbol\lll 136E
 
\newsymbol\ggg 136F
 
\newsymbol\circledS 1073
\newsymbol\pitchfork 1374
\newsymbol\dotplus 1275
\newsymbol\backsim 1376
\newsymbol\backsimeq 1377
\newsymbol\complement 107B
\newsymbol\intercal 127C
\newsymbol\circledcirc 127D
\newsymbol\circledast 127E
\newsymbol\circleddash 127F
\newsymbol\lvertneqq 2300
\newsymbol\gvertneqq 2301
\newsymbol\nleq 2302
\newsymbol\ngeq 2303
\newsymbol\nless 2304
\newsymbol\ngtr 2305
\newsymbol\nprec 2306
\newsymbol\nsucc 2307
\newsymbol\lneqq 2308
\newsymbol\gneqq 2309
\newsymbol\nleqslant 230A
\newsymbol\ngeqslant 230B
\newsymbol\lneq 230C
\newsymbol\gneq 230D
\newsymbol\npreceq 230E
\newsymbol\nsucceq 230F
\newsymbol\precnsim 2310
\newsymbol\succnsim 2311
\newsymbol\lnsim 2312
\newsymbol\gnsim 2313
\newsymbol\nleqq 2314
\newsymbol\ngeqq 2315
\newsymbol\precneqq 2316
\newsymbol\succneqq 2317
\newsymbol\precnapprox 2318
\newsymbol\succnapprox 2319
\newsymbol\lnapprox 231A
\newsymbol\gnapprox 231B
\newsymbol\nsim 231C
\newsymbol\ncong 231D
\newsymbol\diagup 201E
\newsymbol\diagdown 201F
\newsymbol\varsubsetneq 2320
\newsymbol\varsupsetneq 2321
\newsymbol\nsubseteqq 2322
\newsymbol\nsupseteqq 2323
\newsymbol\subsetneqq 2324
\newsymbol\supsetneqq 2325
\newsymbol\varsubsetneqq 2326
\newsymbol\varsupsetneqq 2327
\newsymbol\subsetneq 2328
\newsymbol\supsetneq 2329
\newsymbol\nsubseteq 232A
\newsymbol\nsupseteq 232B
\newsymbol\nparallel 232C
\newsymbol\nmid 232D
\newsymbol\nshortmid 232E
\newsymbol\nshortparallel 232F
\newsymbol\nvdash 2330
\newsymbol\nVdash 2331
\newsymbol\nvDash 2332
\newsymbol\nVDash 2333
\newsymbol\ntrianglerighteq 2334
\newsymbol\ntrianglelefteq 2335
\newsymbol\ntriangleleft 2336
\newsymbol\ntriangleright 2337
\newsymbol\nleftarrow 2338
\newsymbol\nrightarrow 2339
\newsymbol\nLeftarrow 233A
\newsymbol\nRightarrow 233B
\newsymbol\nLeftrightarrow 233C
\newsymbol\nleftrightarrow 233D
\newsymbol\divideontimes 223E
\newsymbol\varnothing 203F
\newsymbol\nexists 2040
\newsymbol\Finv 2060
\newsymbol\Game 2061
\newsymbol\mho 2066
\newsymbol\eth 2067
\newsymbol\eqsim 2368
\newsymbol\beth 2069
\newsymbol\gimel 206A
\newsymbol\daleth 206B
\newsymbol\lessdot 236C
\newsymbol\gtrdot 236D
\newsymbol\ltimes 226E
\newsymbol\rtimes 226F
\newsymbol\shortmid 2370
\newsymbol\shortparallel 2371
\newsymbol\smallsetminus 2272
\newsymbol\thicksim 2373
\newsymbol\thickapprox 2374
\newsymbol\approxeq 2375
\newsymbol\succapprox 2376
\newsymbol\precapprox 2377
\newsymbol\curvearrowleft 2378
\newsymbol\curvearrowright 2379
\newsymbol\digamma 207A
\newsymbol\varkappa 207B
\newsymbol\Bbbk 207C
\newsymbol\hslash 207D
\undefine\hbar
\newsymbol\hbar 207E
\newsymbol\backepsilon 237F

%
\newcount\marknumber	\marknumber=1
\newcount\countdp \newcount\countwd \newcount\countht 
%
%
\ifx\pdfoutput\undefined
\def\rgboo#1{}
\def\postscript#1{\special{" #1}}		
\postscript{
	/bd {bind def} bind def
	/fsd {findfont exch scalefont def} bd
	/sms {setfont moveto show} bd
	/ms {moveto show} bd
	/pdfmark where		
	{pop} {userdict /pdfmark /cleartomark load put} ifelse
	[ /PageMode /UseOutlines		
	/DOCVIEW pdfmark}
\def\bookmark#1#2{\postscript{		
	[ /Dest /MyDest\the\marknumber /View [ /XYZ null null null ] /DEST pdfmark
	[ /Title (#2) /Count #1 /Dest /MyDest\the\marknumber /OUT pdfmark}%
	\advance\marknumber by1}
\def\pdfclink#1#2#3{%
	\hskip-.25em\setbox0=\hbox{#2}%
		\countdp=\dp0 \countwd=\wd0 \countht=\ht0%
		\divide\countdp by65536 \divide\countwd by65536%
			\divide\countht by65536%
		\advance\countdp by1 \advance\countwd by1%
			\advance\countht by1%
		\def\linkdp{\the\countdp} \def\linkwd{\the\countwd}%
			\def\linkht{\the\countht}%
	\postscript{
		[ /Rect [ -1.5 -\linkdp.0 0\linkwd.0 0\linkht.5 ] 
		/Border [ 0 0 0 ]
		/Action << /Subtype /URI /URI (#3) >>
		/Subtype /Link
		/ANN pdfmark}{\rgb{#1}{#2}}}
%
%
\else
\def\rgboo#1{\pdfliteral{#1 rg #1 RG}}
\pdfcatalog{/PageMode /UseOutlines}		
\def\bookmark#1#2{
	\pdfdest num \marknumber xyz
	\pdfoutline goto num \marknumber count #1 {#2}
	\advance\marknumber by1}
\def\pdfklink#1#2{%
	\noindent\pdfstartlink user
		{/Subtype /Link
		/Border [ 0 0 0 ]
		/A << /S /URI /URI (#2) >>}{\rgb{1 0 0}{#1}}%
	\pdfendlink}
\fi

\def\rgbo#1#2{\rgboo{#1}#2\rgboo{0 0 0}}
\def\rgb#1#2{\mark{#1}\rgbo{#1}{#2}\mark{0 0 0}}
\def\pdfklink#1#2{\pdfclink{1 0 0}{#1}{#2}}
\def\pdflink#1{\pdfklink{#1}{#1}}
%
%
\newcount\seccount  
\newcount\subcount  
\newcount\clmcount  
\newcount\equcount  
\newcount\refcount  
\newcount\demcount  
\newcount\execount  
\newcount\procount  
\seccount=0
\equcount=1
\clmcount=1
\subcount=1
\refcount=1
\demcount=0
\execount=0
\procount=0
%
\def\proof{\medskip\noindent{\bf Proof.\ }}
\def\proofof(#1){\medskip\noindent{\bf Proof of \csname c#1\endcsname.\ }}
\def\qed{\hfill{\sevenbf QED}\par\medskip}
\def\references{\bigskip\noindent\hbox{\bf References}\medskip
                \ifx\pdflink\undefined\else\bookmark{0}{References}\fi}
\def\addref#1{\expandafter\xdef\csname r#1\endcsname{\number\refcount}
    \global\advance\refcount by 1}

\def\nextremark #1\par{\item{$\circ$} #1}
\def\firstremark #1\par{\bigskip\noindent{\bf Remarks.}
     \smallskip\nextremark #1\par}
\def\abstract#1\par{{\baselineskip=10pt
    \eightpoint\narrower\noindent{\eightbf Abstract.} #1\par}}
%
\def\equtag#1{\expandafter\xdef\csname e#1\endcsname{(\number\seccount.\number\equcount)}
              \global\advance\equcount by 1}
\def\equation(#1){\equtag{#1}\eqno\csname e#1\endcsname}
\def\equ(#1){\hskip-0.03em\csname e#1\endcsname}
%
\def\clmtag#1#2{\expandafter\xdef\csname cn#2\endcsname{\number\seccount.\number\clmcount}
                \expandafter\xdef\csname c#2\endcsname{#1~\number\seccount.\number\clmcount}
                \global\advance\clmcount by 1}
\def\claim #1(#2) #3\par{\clmtag{#1}{#2}
    \vskip.1in\medbreak\noindent
    {\bf \csname c#2\endcsname .\ }{\sl #3}\par
    \ifdim\lastskip<\medskipamount
    \removelastskip\penalty55\medskip\fi}
\def\clm(#1){\csname c#1\endcsname}
\def\clmno(#1){\csname cn#1\endcsname}
%
\def\sectag#1{\global\advance\seccount by 1
              \expandafter\xdef\csname sectionname\endcsname{\number\seccount. #1}
              \equcount=1 \clmcount=1 \subcount=1 \execount=0 \procount=0}
\def\section#1\par{\vskip0pt plus.1\vsize\penalty-40
    \vskip0pt plus -.1\vsize\bigskip\bigskip
    \sectag{#1}
    \message{\sectionname}\leftline{\magtenbf\sectionname}
    \nobreak\smallskip\noindent
    \ifx\pdflink\undefined
    \else
      \bookmark{0}{\sectionname}
    \fi}
%
\def\subtag#1{\expandafter\xdef\csname subsectionname\endcsname{\number\seccount.\number\subcount. #1}
              \global\advance\subcount by 1}
\def\subsection#1\par{\vskip0pt plus.05\vsize\penalty-20
    \vskip0pt plus -.05\vsize\medskip\medskip
    \subtag{#1}
    \message{\subsectionname}\leftline{\tenbf\subsectionname}
    \nobreak\smallskip\noindent
    \ifx\pdflink\undefined
    \else
      \bookmark{0}{.... \subsectionname}  
    \fi}
%
\def\demtag#1#2{\global\advance\demcount by 1
              \expandafter\xdef\csname de#2\endcsname{#1~\number\demcount}}
\def\demo #1(#2) #3\par{
  \demtag{#1}{#2}
  \vskip.1in\medbreak\noindent
  {\bf #1 \number\demcount.\enspace}
  {\rm #3}\par
  \ifdim\lastskip<\medskipamount
  \removelastskip\penalty55\medskip\fi}
\def\dem(#1){\csname de#1\endcsname}
%
\def\exetag#1{\global\advance\execount by 1
              \expandafter\xdef\csname ex#1\endcsname{Exercise~\number\seccount.\number\execount}}
\def\exercise(#1) #2\par{
  \exetag{#1}
  \vskip.1in\medbreak\noindent
  {\bf Exercise \number\execount.}
  {\rm #2}\par
  \ifdim\lastskip<\medskipamount
  \removelastskip\penalty55\medskip\fi}
\def\exe(#1){\csname ex#1\endcsname}
%
\def\protag#1{\global\advance\procount by 1
              \expandafter\xdef\csname pr#1\endcsname{\number\seccount.\number\procount}}
\def\problem(#1) #2\par{
  \ifnum\procount=0
    \parskip=6pt
    \vbox{\bigskip\centerline{\bf Problems \number\seccount}\nobreak\medskip}
  \fi
  \protag{#1}
  \item{\number\procount.} #2}
\def\pro(#1){Problem \csname pr#1\endcsname}
%
%
%
\def\rightheadline{\hfil}
\def\leftheadline{\sevenrm\hfil HANS KOCH\hfil}
\headline={\ifnum\pageno=\firstpage\hfil\else
\ifodd\pageno{{\fiverm\rightheadline}\number\pageno}
\else{\number\pageno\fiverm\leftheadline}\fi\fi}
\footline={\ifnum\pageno=\firstpage\hss\tenrm\folio\hss\else\hss\fi}

\let\cl=\centerline

\let\eps=\varepsilon
\let\sss=\scriptscriptstyle

\def\AA{{\cal A}}
\def\BB{{\cal B}}

\def\HH{{\cal H}}
\def\II{{\cal I}}

\def\KK{{\cal K}}
\def\LL{{\cal L}}

\def\NN{{\cal N}}

\def\PP{{\cal P}}

\def\RR{{\cal R}}
\def\SS{{\cal S}}

\def\VV{{\cal V}}

\def\ssN{{\sss N}}

\def\rmC{{\rm C}}
\def\rmL{{\rm L}}
\def\id{{\rm I}}

%
\newfam\dsfam
\def\mathds #1{{\fam\dsfam\tends #1}}

\font\tends=dsrom10
\font\eightds=dsrom8
\textfont\dsfam=\tends
\scriptfont\dsfam=\eightds
%

\def\integer{{\mathds Z}}

\def\real{{\mathds R}}
\def\complex{{\mathds C}}

\def\proj{{\Bbb P}}

\def\Id{{\Bbb I}}

\def\defeq{\mathrel{\mathop=^{\sss\rm def}}}
\def\half{{1\over 2}}

\def\quarter{{1\over 4}}

\def\tquarter{{\textstyle\quarter}}
%

%

%

%

%

\input miniltx

\ifx\pdfoutput\undefined
  \def\Gin@driver{dvips.def}  
\else
  \def\Gin@driver{pdftex.def} 
\fi
 
\input graphicx.sty
\resetatcatcode
%
\font\twelvelfb=lfb12
\def\Chi{\hbox{\twelvelfb c}}
\def\buN{{\hbox{\magnineeufm N}}}
\def\ssN{{\sss N}}
\let\disk\Omega
\def\sphere{{\Bbb S}}
\def\SO{{\hbox{\rm SO}}}
\def\SU{{\hbox{\rm SU}}}
\def\rmH{{\rm H}}
\def\evensym{0}
\def\oddsym{1}
\def\paritysym{\sigma}
\def\even{{\hskip1pt\evensym}}
\def\odd{{\hskip1pt\oddsym}}
\def\parity{{\hskip1pt\paritysym}}
\def\range{\hskip1pt{-}\hskip0.6pt}

\def\Langle{\bigl\langle}
\def\Rangle{\bigr\rangle}
\def\uzero{\bar u}
\def\boxit#1#2{\leavevmode\hbox{\vrule
\vtop{\vbox{\hrule\kern#1\hbox{\kern#1#2\kern#1}}\kern#1\hrule}\vrule}}
\addref{CNZ}
\addref{SSW}
\addref{PW}
\addref{BaCa}
\addref{Koo}
\addref{Kin}
\addref{Tango}
\addref{Kanjin}
\addref{Wue}
\addref{Jan}
\addref{AKi}
\addref{AKii}
\addref{CZ}
\addref{WaNa}
\addref{CLJii}
\addref{FH}
\addref{FL}
\addref{BVCDLVWY}
\addref{AKiii}
\addref{More}
\addref{Regge}
\addref{RBMW}
\addref{RaYu}
\addref{JoFo}
\addref{Files}
\addref{Ada}
\addref{Gnat}
\addref{IEEE}
\addref{MPFR}
%
\def\leftheadline{\sevenrm\hfil
GIANNI ARIOLI and HANS KOCH\hfil}
\def\rightheadline{\sevenrm\hfil semilinear equations on the disk\hfil}
\cl{{\magtenbf Validated numerical solutions}}
\cl{{\magtenbf for some semilinear elliptic equations on the disk}}
\bigskip

\cl{Gianni Arioli
\footnote{$^1$}
{{\sevenrm Department of Mathematics, Politecnico di Milano,
Piazza Leonardo da Vinci 32, 20133 Milano.}}
and Hans Koch
\footnote{$^2$}
{{\sevenrm Department of Mathematics, University of Texas at Austin,
Austin, TX 78712}}
}

\bigskip
\abstract
Starting with approximate solutions
of the equation $-\Delta u=wu^3$ on the disk,
with zero boundary conditions,
we prove that there exist true solutions nearby.
One of the challenges here lies in the fact
that we need simultaneous and accurate control
of both the (inverse) Dirichlet Laplacean and nonlinearities.
We achieve this with the aid of a computer,
using a Banach algebra of real analytic functions,
based on Zernike polynomials.
Besides proving existence, and symmetry properties,
we also determine the Morse index of the solutions.

\section Introduction

In this paper we consider semilinear elliptic equations of the form
$$
-\Delta u=wf'(u)\,,\qquad u\bigm|_{\partial\disk}\,=0\,,
\equation(Main)
$$
where $\disk$ is the unit disk in $\real^2$,
$w$ is a nonnegative function on $\disk$,
and $f'$ is the derivative of a regular function $f$ on $\real$.
Our primary goal is to develop techniques
that can be used to prove the existence of solutions
in a constructive way, with the help of a computer.
In the concrete cases considered here,
$w$ is always radial (invariant under rotations) and $f'(u)=u^3$.
But it will be clear from our description
that the same methods work for other choices of $w$ and $f$.
In fact, similar techniques should apply to other types of equations,
and to other radially symmetric domains in $\real^2$ and $\real^3$.

Before giving more details,
let us state a result that will help to set the stage.

\claim Theorem(HillRing)
There exists a positive radial polynomial $w$ on $\disk$,
such that the equation \equ(Main) with $f'(u)=u^3$
admits a real analytic solution $u=u_w$ that has Morse index $2$,
with the property that $|u_w|$ is not invariant
under any nontrivial rotation.

The weight function $w$ and the solution $u_w$ are shown in Figure 1.
A precise definition of $w$ is given in [\rFiles].
We note that $u_w$ is symmetric under a reflection.
This is one symmetry that solutions cannot avoid [\rPW].
Our goal was to find an index-$2$ solution that has no other symmetries.

Concerning the Morse index,
recall that solutions of equation \equ(Main)
are critical points of the functional $J$ on $\rmH^1_0(\disk)$,
$$
J(u)=\int_\disk\,\Bigl[{\textstyle\half}
\bigl|\nabla u\bigr|^2-wf(u)\Bigr]\,dxdy\,,
\equation(Ju)
$$
assuming that $f$ satisfies some growth and regularity conditions.
The Morse index of a critical point $u$
is the number of descending directions of $J$ at $u$.

One of the difficulties with proving \clm(HillRing) is that $\disk$ is a disk.
For a square domain, an analogous result was proved in [\rAKi].
And for the disk, it is possible [\rAKii] to obtain an accurate numerical ``solution''
that looks as shown in Figure 1.
But we have hitherto been unable to prove that there exists a true solution nearby.

\vskip0.5cm
\hbox{\hskip0.7cm
\includegraphics[height=5.0cm,width=6.2cm]{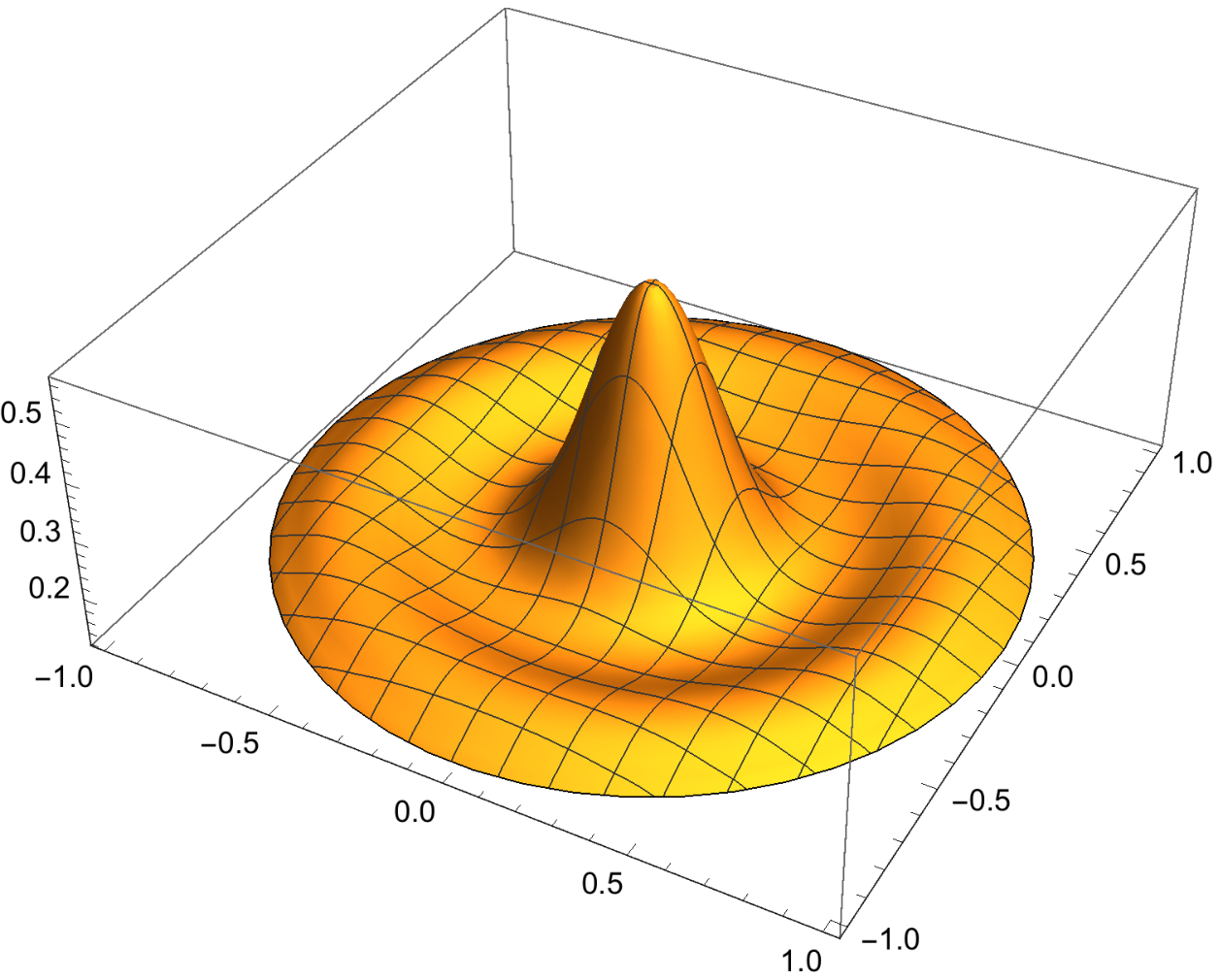}
\includegraphics[height=5.0cm,width=6.2cm]{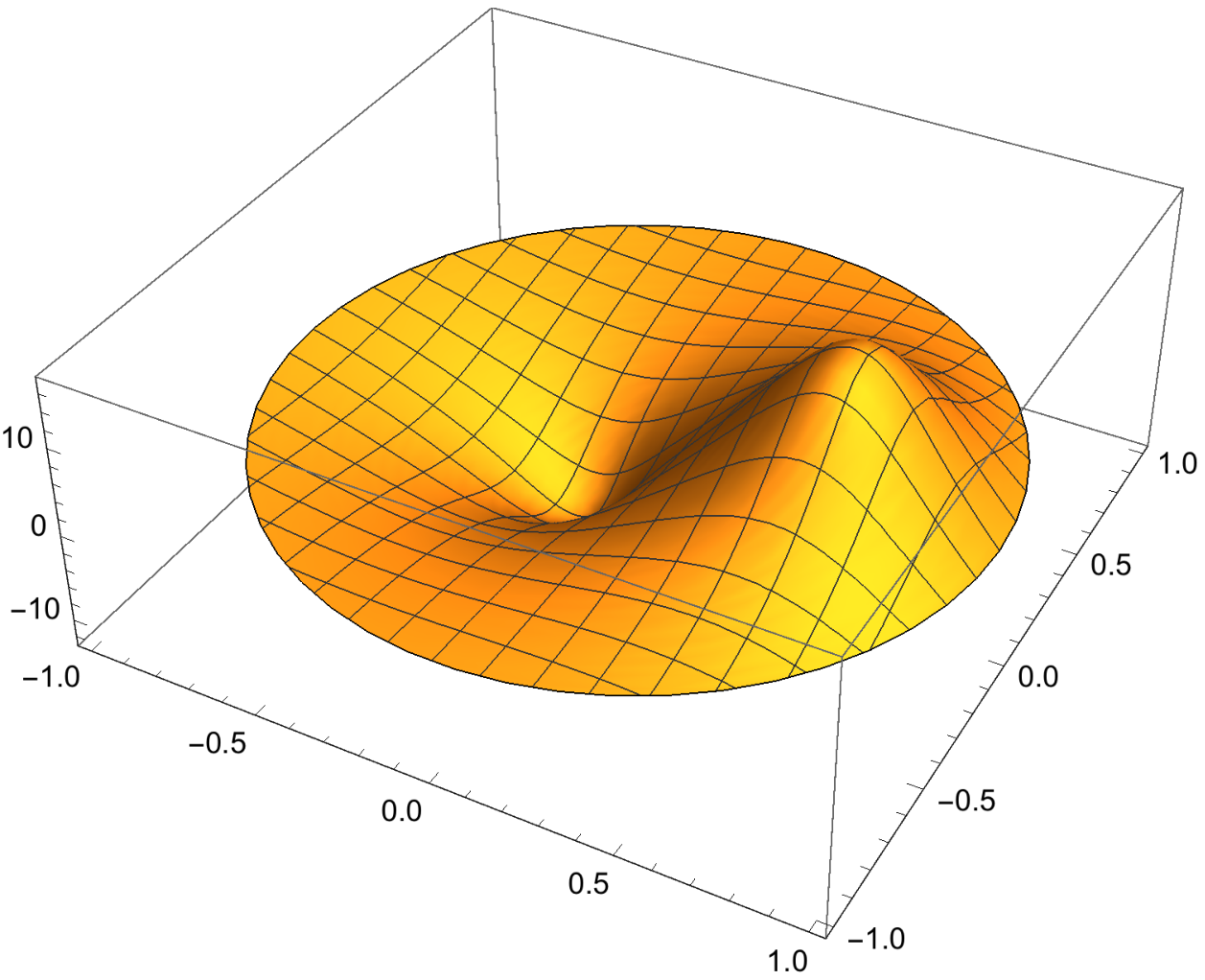}}
\centerline{\eightpoint{\noindent\bf Figure 1.}
The weight function $w$ and solution $u_w$ described in \clm(HillRing).}
\vskip0.5cm

Before describing our approach in more detail,
let us state two other results that can be proved in a similar way.
The first results concern again
``minimally symmetric solutions to a highly symmetric problem''.
While the weight $w$ in \clm(HillRing) had to be chosen carefully
to obtain a minimally symmetric solution of index $2$,
a standard H\'enon weight $w(r,\vartheta)=r^\alpha$
suffices in the index-$1$ case.
Here, and in what follows, $(r,\vartheta)$ denote the standard
polar coordinates on $\disk$.

\claim Theorem(HenonOne)
For $\alpha=2,4,6$, the equation \equ(Main), with $w=r^\alpha$ and $f'(u)=u^3$,
admits a real analytic solution $u=u_\alpha>0$. This solution has Morse index $1$
and is not invariant under any nontrivial rotation.

The solutions $u_2$, $u_4$, and $u_6$ are shown in Figure 2.

\vskip0.5cm
\hbox{\hskip0.7cm
\includegraphics[height=3.5cm,width=4cm]{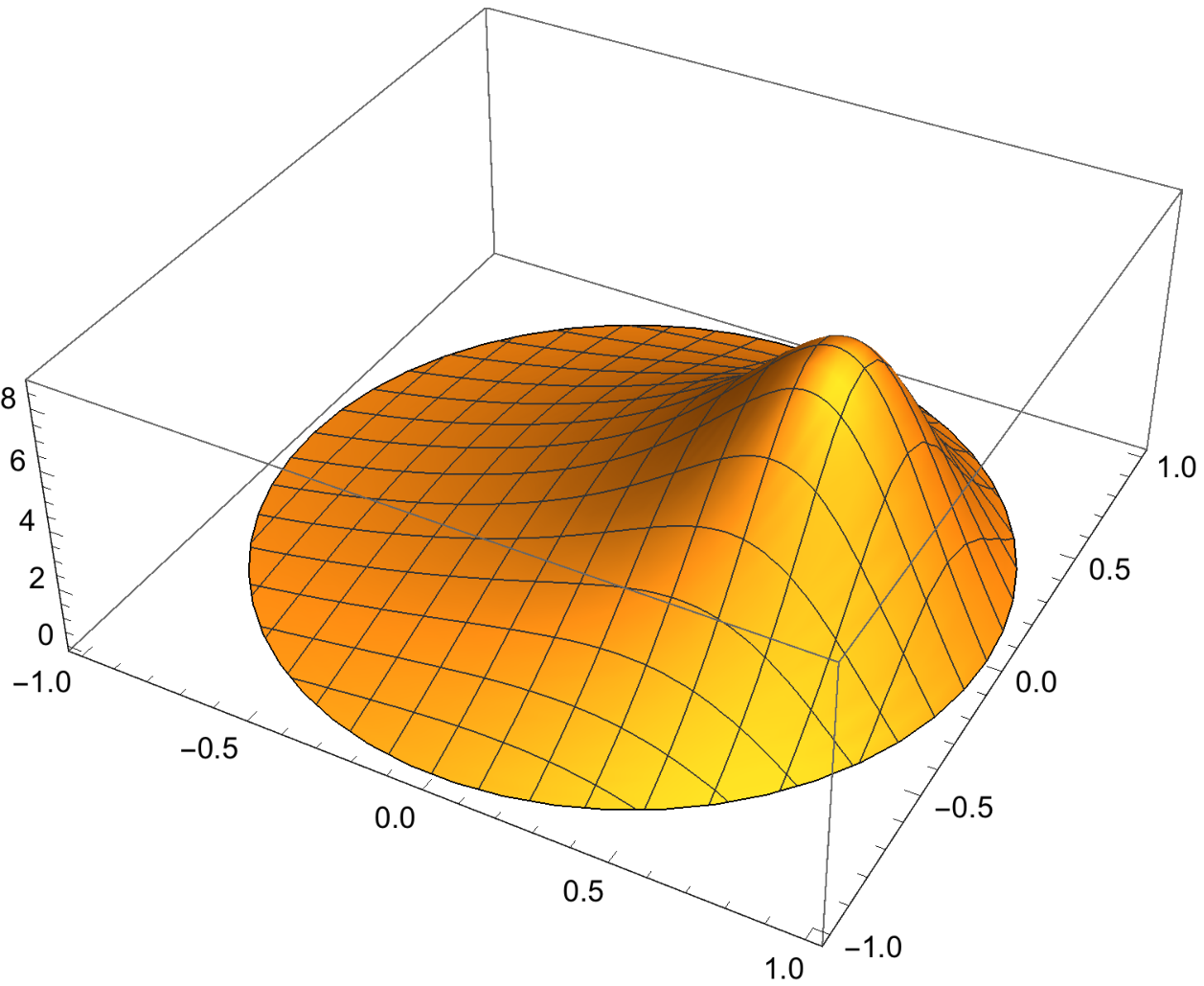}
\includegraphics[height=3.5cm,width=4cm]{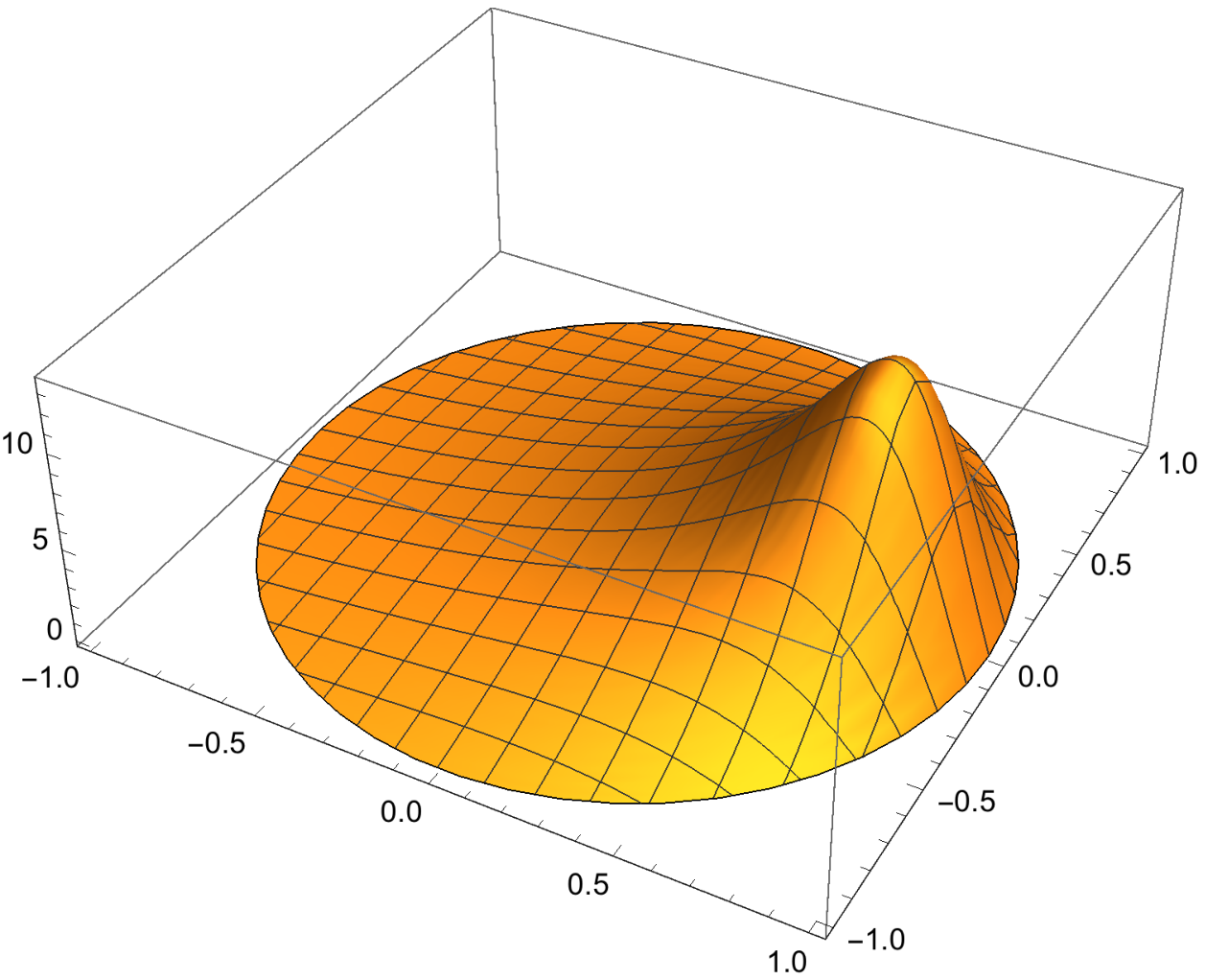}
\includegraphics[height=3.5cm,width=4cm]{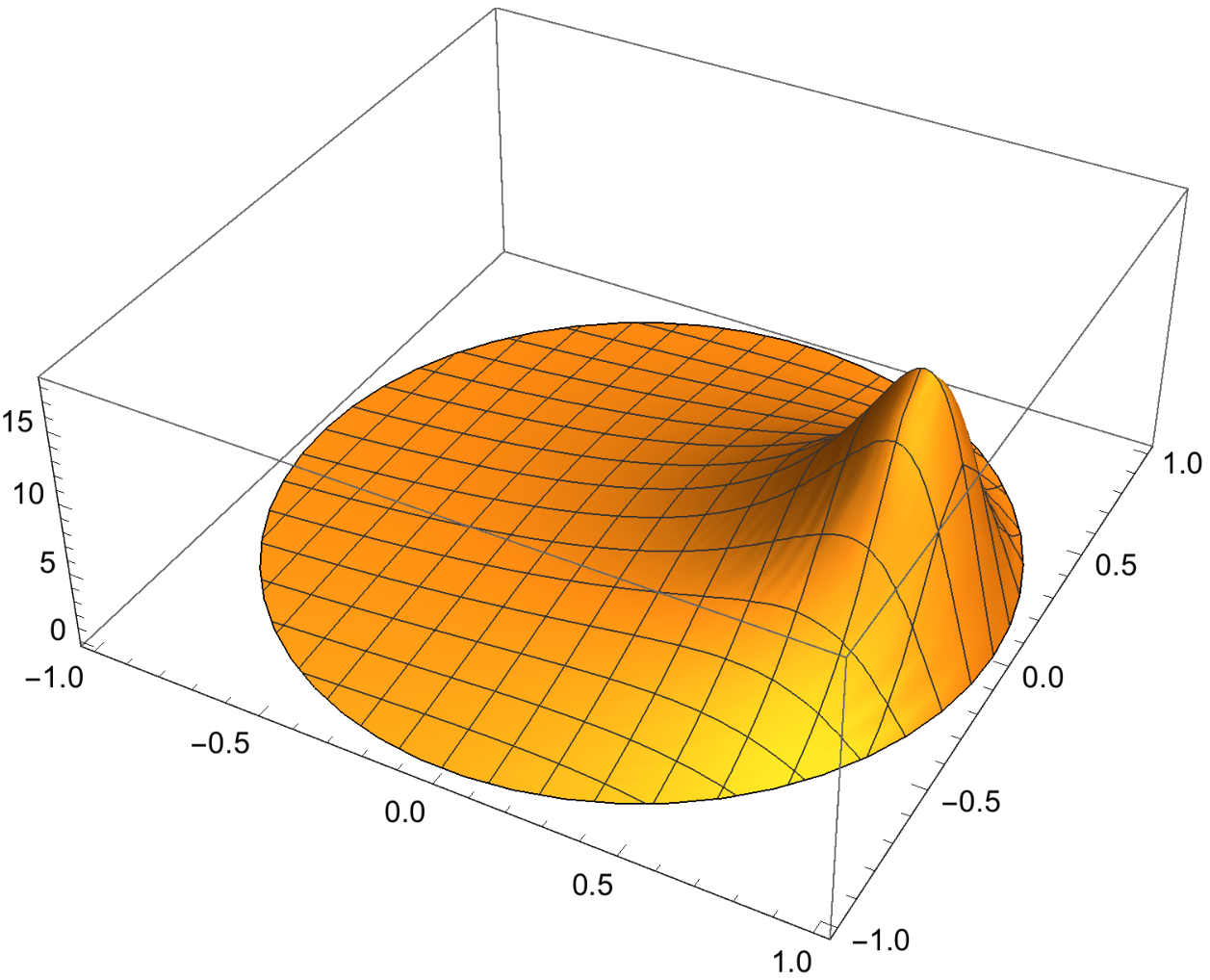}}
\centerline{\eightpoint{\noindent\bf Figure 2.}
The solutions $u_2$, $u_4$, and $u_6$ described in \clm(HenonOne).}
\vskip0.5cm

The same result, but without the statement about the lack of symmetry,
is easy to prove: minimizing $J$ on the Nehari manifold
$\buN=\bigl\{u\in\rmH^1_0(\disk):\,DJ(u)u=0\,,u\ne 0\bigr\}$
shows that index-$1$ solutions exist
and that they do not vanish anywhere on $\disk$.
Intuitively, the asymmetry of the positive minimizers $u_\alpha$
stems from fact that the term $-wf(u)$ in the integral \equ(Ju)
rewards $u$ for concentrating off-center.
Indeed, an analogue of \clm(HenonOne) can be proved by variational methods
for sufficiently large values of $\alpha$ [\rSSW];
see also [\rBaCa] and references therein.
Numerical results on a number of nonlinear elliptic equation
can be found in [\rCNZ].
They include positive non-radial solutions $u_\alpha$ as described in 
\clm(HenonOne), but for $\alpha=1,9$.
Given that the positive solution for $\alpha=0$ is radial,
one expects that there is a symmetry-breaking bifurcation
as $\alpha$ is increased from $0$ to $1$.

Our method of proof is not limited to solutions of index $1$ or $2$,
although the computations become impractical at high index.
In the next theorem,
we consider two solutions that are close to sums of index-$1$ solutions,
$$
u_{\alpha,n}\approx\sum_{m=1}^{2n}(S_n)^m u_\alpha\,,\qquad
(S_nu)(r,\vartheta)=-u(r,\vartheta+\pi/n)\,.
\equation(ualphan)
$$
If $u_\alpha$ is one of the solutions described in \clm(HenonOne),
then the functions in the above sum are solutions of the same equation;
and if $n$ is not too large,
then most of their mass is contained in mutually disjoint sectors of the disk.
Thus, the sum in \equ(ualphan)
is an approximate solution of the H\'enon equation,
and we expect to find a true solution nearby.
Furthermore, this solution should have index $2n$.
Indeed, this holds in the two cases considered here:

\claim Theorem(HenonTwo)
For $n=1,2$, the equation \equ(Main), with $w=r^{2}$ and $f'(u)=u^3$,
admits a non-radial real analytic solution $u=u_{2,n}$
that is invariant under $S_n$ and has index $2n$.

The functions $u_{2,1}$ and $u_{2,2}$ are shown in Figure 3.
We expect that solutions of the type \equ(ualphan) exist for any given $n>0$,
provided that $\alpha$ is sufficiently large.

\vskip0.3cm
\hbox{\hskip0.7cm
\includegraphics[height=5.0cm,width=6.2cm]{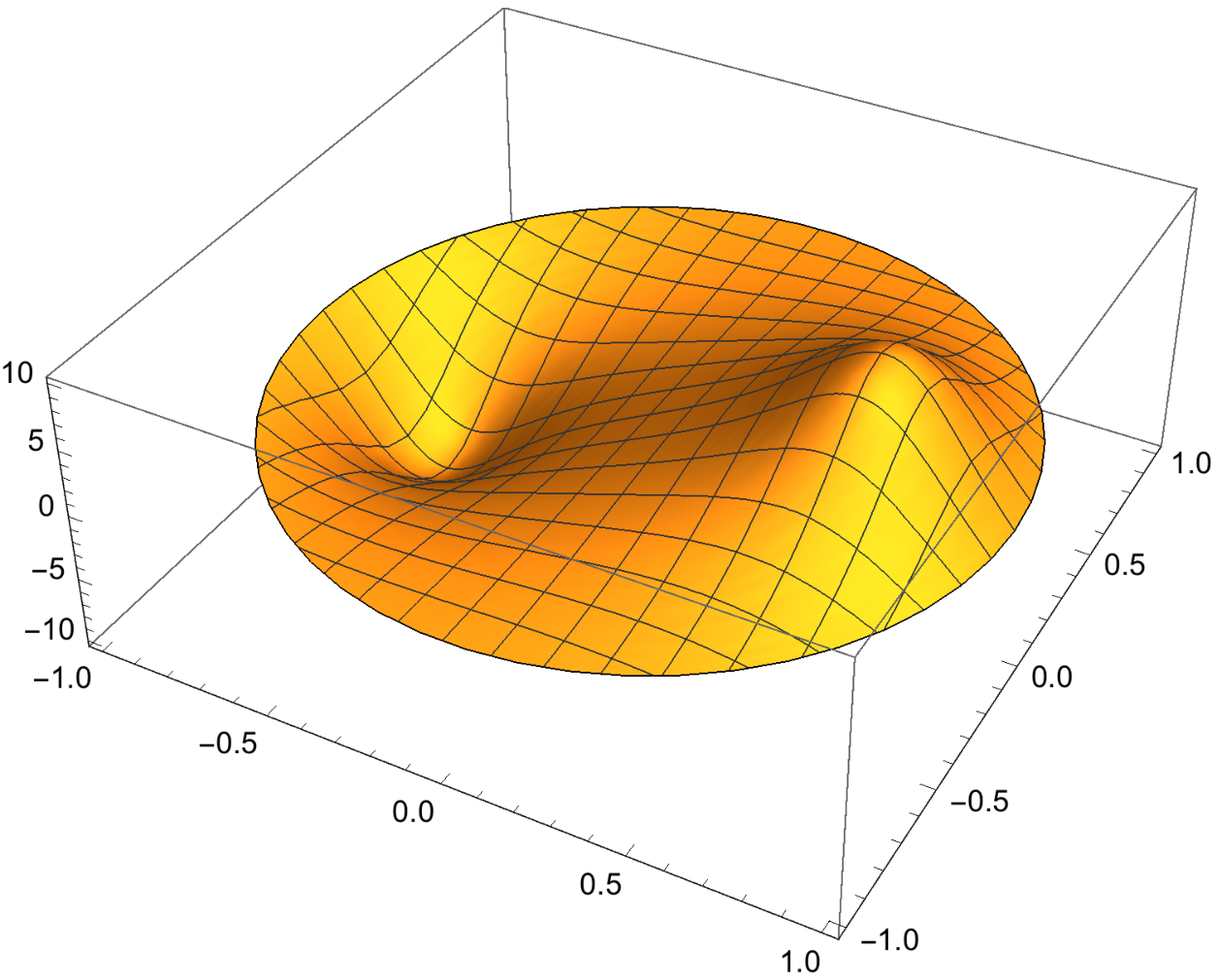}
\includegraphics[height=5.0cm,width=6.2cm]{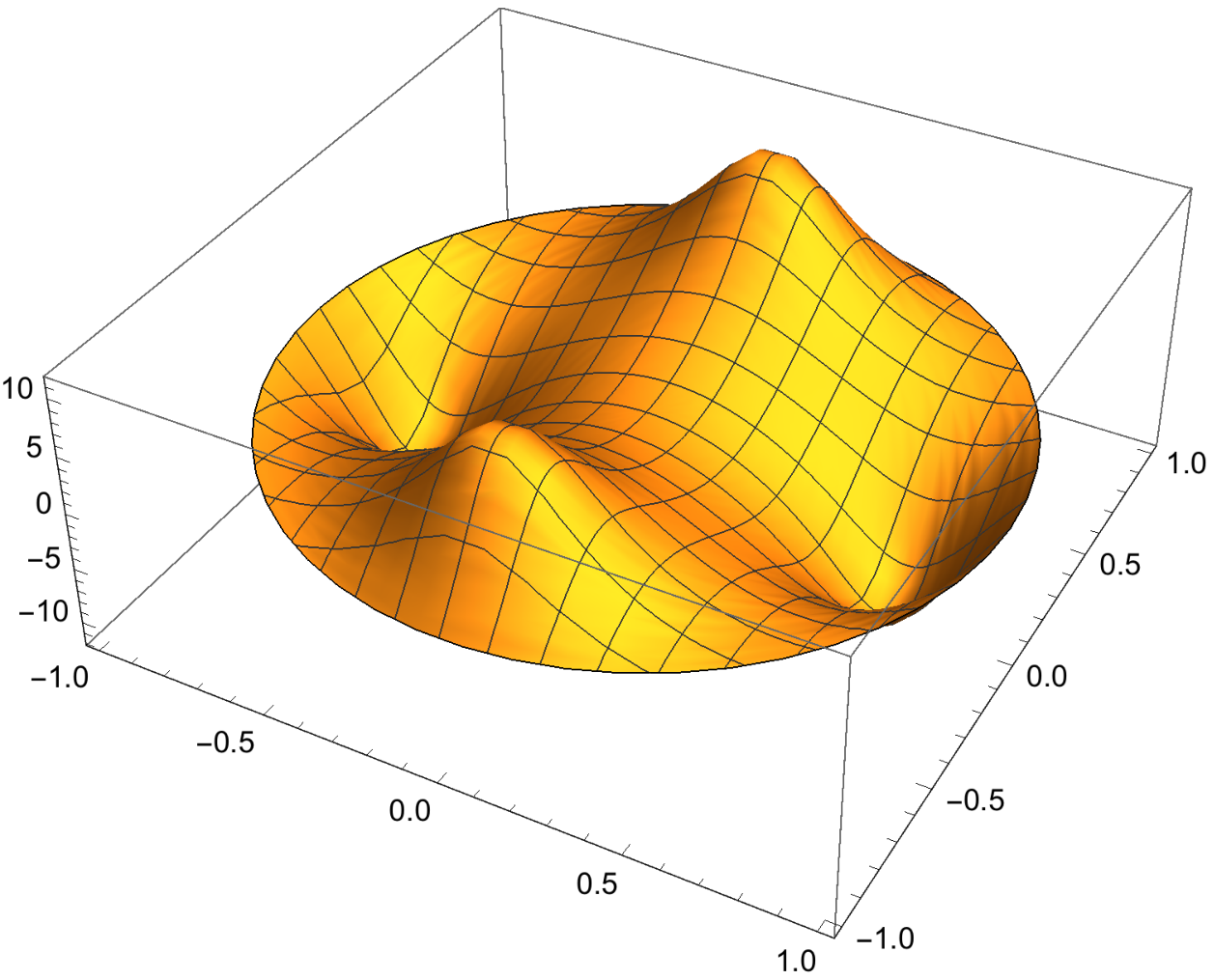}}
\vskip0.3cm

\centerline{\eightpoint{\noindent\bf Figure 3.}
The solutions $u_{2,1}$ and $u_{2,2}$ described in \clm(HenonTwo).}
\vskip0.4cm

As indicated earlier, the three theorems stated above
are proved with the aid of a computer.
In addition to the properties described in these theorems,
we obtain accurate bounds on the difference $u-\uzero$
between the true solution $u$ and a numerical approximation $\uzero$.
The accuracy of the result is limited only by the computational resources available.
To give a rough idea: our estimates on the solution $u=u_w$
described in \clm(HillRing), carried out on a standard desktop machine,
yield an upper bound less than $2^{-35}$ on the norm (defined later) of $u-\uzero$
relative to the norm of $\uzero$.

Following a strategy that has been successful
in many other computer-assisted proofs in analysis [\rAKi\range\rMore],
we start by converting the given equation \equ(Main)
to a fixed point equation for a suitable operator $G$.
As in [\rAKi], we use
$$
G(u)=-\Delta^{-1}\bigl[wf'(u)\bigr]\,,
\equation(MainFix)
$$
where $\Delta^{-1}$ is the inverse Dirichlet Laplacean on $\Omega$.
Then we consider a Newton-type map $\NN$
associated with $G$ and prove that $\NN$ is a contraction
in a small neighborhood of an approximate fixed point $\uzero$.

Clearly, approximations play a crucial role.
Without loss of generality,
we are looking for a representation $u=\sum_k c_k\Psi_k$,
where approximation corresponds to truncating the series to a finite sum.
So it is desirable to work with a space $\BB$
and basis function $\Psi_k$ that

\item{$(a)$} {\sl are well adapted to the operators involved},
\item{$(b)$} {\sl have useful algebraic properties and}
\item{$(c)$} {\sl good approximation properties}.

\smallskip\noindent
The same criteria apply to most computer-assisted proofs in analysis
[\rAKi\range\rMore].
In problems that involve a Laplacean with a compact inverse,
the ideal way to satisfy $(a)$ is to take for $\Psi_k$
the $k$-th eigenvector of $\Delta^{-1}$.
This works well for rectangular domains and Fourier series.
In this case, the basis functions $\Psi_k$
also have a simple product expansion $\Psi_i\Psi_j=\sum_k c_{i,j,k}\Psi_k$.
This is a desired property $(b)$
in problems that include nonlinearities such as the term
$f'(u)$ in \equ(MainFix).
Ideally, $\BB$ is a Banach algebra.
Concerning $(c)$, the expected solution $u$
should have coefficients $c_k$
that decrease more rapidly than those for a typical function in $\BB$.

In the problems considered here,
the eigenfunctions of the Dirichlet Laplacean (on the disk) are of the form
$\Psi_k(r,\vartheta)=\psi_k(r)e^{im_k\vartheta}$,
where $\psi_k$ is an appropriately scaled Bessel function.
Unfortunately, there is no convenient product expansion for these Bessel functions.
Thus, in earlier attempts to prove a result like \clm(HillRing),
we used various approximations or alternatives for the Bessel functions.
Some choices worked well for numerical computations,
as described in [\rAKii], but they never led to a successful proof.
We had also considered using Zernike polynomials,
but obviously not carefully enough.

As it turns out, the Zernike polynomials are
close to ideal for the type of problems considered here.
There are good reasons for this, as we will explain below.
The Zernike polynomials $R^m_n$ are widely used in optics.
But despite the vast literature on (the use of) these polynomials,
we found no evidence that would have justified trying yet another approach.
That is, until we became aware of the references [\rTango,\rJan].

In [\rTango] it is shown that the Zernike functions
$V^m_n(r,\vartheta)=R^m_n(r)e^{im\vartheta}$ have a product expansion
whose coefficients are the squares of certain Clebsch-Gordan coefficients.
This property is obtained by relating 
the Zernike functions to generalized spherical harmonics,
using the azimuthal projection of the sphere to the disk.
So in essence, one works indirectly with functions on $\SO(3)$,
whence the nice behavior under multiplication (product representations).
But the Laplacean on the sphere does not map
to the Laplacean on the disk, which creates a potential conflict
between the desired properties $(a)$ and $(b)$.
Surprisingly, this problem is quite harmless:
both $\Delta V^m_n$ and $\Delta^{-1}V^m_n$
are linear combinations of at most three Zernike functions [\rJan].

This motivates the following expansion
for our functions on the unit disk:
$$
u(r,\vartheta)=\sum_{m,l=0}^\infty
R^m_{m+2l}(r)\Bigl[a_{m,l}\cos(m\vartheta)+b_{m,l}\sin(m\vartheta)\Bigr]\,,
\equation(uZernikeSeries)
$$
with $b_{0,l}=0$ for all $l$.
To be more precise, the solutions of \equ(Main) described in the preceding theorems
are symmetric under the reflection $\vartheta\mapsto-\vartheta$,
so their coefficients $b_{m,l}$ all vanish.
The function spaces used in our analysis are the following.
Given $\rho>1$, we define $\BB_\rho$
to be the space of all functions \equ(uZernikeSeries)
that have a finite norm
$$
\|u\|_\rho=\sum_{m,l=0}^\infty\bigl(|a_{m,l}|+|b_{m,l}|\bigr)\rho^{m+2l}\,.
\equation(BBrhoNormCS)
$$
Using the relationship between the Zernike functions
and the generalized spherical harmonics,
one immediately gets complex bounds on the Zernike functions.
This, together with the Clebsch-Gordan series
for products can be used to show that $\BB_\rho$
is in fact a Banach algebra of real analytic functions on $\disk$.

Clearly, if $u\in\BB_\varrho$ with $\varrho>\rho$,
then $u$ is well approximated in $\BB_\rho$ by truncated sums \equ(uZernikeSeries).
This allows us to obtain highly accurate numerical approximations to the solutions
described in Theorems \clmno(HillRing), \clmno(HenonOne), and \clmno(HenonTwo).
The limiting factor here is the computation or storage of Clebsch-Gordan coefficients.
The computation is quite costly,
so we compute the coefficients beforehand and store them in an array.
To reduce symmetries and thus storage space, we use ideas described in [\rRaYu].

\section Zernike functions

Here we introduce the Zernike functions
and describe the properties that we need in our analysis.
Our need for the product expansion and complex bounds
favors the following approach.

Consider unitary representations of $\SU(2)$,
with Hermitian generators $L_1$, $L_2$, and $L_3$
satisfying $[L_2,L_3]=iL_1$, $[L_3,L_1]=iL_2$, and $[L_1,L_2]=iL_3$.
The eigenvalues of $L_3$ are commonly referred to as weights.
Each irreducible representation is characterized
uniquely (up to unitary equivalence) by the value $\nu$ of the largest weight,
which is a nonnegative half-integer.
In such a representation, the spectrum of each generator $L_j$
is the set $W_\nu=\{-\nu,-\nu+1,\ldots,\nu-1,\nu\}$,
and all the eigenvalues are simple.
Furthermore, the $d=2\nu+1$ eigenvectors of $L_3$
constitute an orthogonal basis for the underlying space.
Using the bra-ket notation that is common in physics,
the normalized eigenvector of $L_3$ with eigenvalue $\mu\in W_\nu$
will be denoted by $|\nu\mu\rangle$.

Without loss of generality,
we may identify each operator $R$ in our representation
by the $d\times d$ matrix 
whose elements are $R_{\mu,\mu'}=\langle\nu\mu|R|\nu\mu'\rangle$.
Here we consider the matrices $e^{-i\alpha L_3}e^{-i\beta L_2}e^{-i\gamma L_3}$,
also known as the Wigner $D$-matrices.
The angles $(\alpha,\beta,\gamma)$  describe the orientation
of a coordinate frame in $\real^3$ with respect to a fixed frame.
We will restrict to $\gamma=0$,
which suffices to describe the direction of the rotated $3$-axis.
Then the matrix elements
$$
D^\nu_{\mu,\mu'}(\alpha,\beta)\defeq
\Langle\nu\mu\bigl|e^{-i\alpha L_3}e^{-i\beta L_2}\bigr|\nu\mu'\Rangle
=e^{-i\mu\alpha}\Langle\nu\mu\bigl|e^{-i\beta L_2}\bigr|\nu\mu'\Rangle
\equation(Dnualphabeta)
$$
define functions on the unit sphere,
also referred to as generalized spherical harmonics.
The functions $D^\nu_{\mu,0}$ are in essence the ordinary spherical harmonics.

Consider now the diagonal elements $D^\nu_{\mu,\mu}$.
They are even functions of $\beta$.
This follows e.g.~from the identity
$(-1)^{\nu-L_3}L_2(-1)^{\nu-L_3}=-L_2$, which can be derived by elementary computations.
Using the spectral decomposition $L_2=\sum_\kappa\kappa P_\kappa$
of the operator $L_2$, we may write
$$
D^\nu_{\mu,\mu}(\alpha,\beta)
=e^{-i\mu\alpha} \sum_{\kappa\in W_\nu}
\langle\nu\mu|P_\kappa|\nu\mu\rangle\cos(\kappa\beta)\,.
\equation(WignerDef)
$$
Notice that the sum in this equation is a polynomial
of degree $2\nu$ in the variable $\cos(\beta/2)$.
Let $n=2\nu$ and $m=2\mu$.
Then we have
$$
V^m_n(r,\vartheta)\defeq
D^\nu_{\mu,\mu}(-2\vartheta,\beta)
=R^m_n(r)e^{im\vartheta}\,,\qquad r=\cos(\beta/2)\,,
\equation(VRfromD)
$$
where $R^m_n$ is a polynomial of degree $n$.
This equation will serve as our definition
of the Zernike functions $V^m_n$ and the Zernike polynomials $R^m_n$.

\smallskip
Next, we describe a few properties of the functions $V^m_n$
that will be needed in our analysis.
Using the well-known Clebsch-Gordan series for the product
of Wigner functions $D^\nu_{\mu,\mu'}$,
one obtains directly the product rule
$$
V_{n_1}^{m_1}V_{n_2}^{m_2}
=\sum_{n_3}\bigl|\langle\nu_1\mu_1\nu_2\mu_2|\nu_3\mu_3\rangle\bigr|^2
\,V_{n_3}^{m_3}\,,\qquad m_3=m_1+m_2\,.
\equation(VProdRule)
$$
Here, $\langle\nu_1\mu_1\nu_2\mu_2|\nu_3\mu_3\rangle$
are the so-called Clebsch-Gordan coefficients,
with $\nu_j=m_j/2$ and $\mu_j=m_j/2$.
These coefficients will be discussed in Section 5.

Another immediate consequence of the definition \equ(VRfromD)
is the following complex bound.
Given that $-\nu\le L_2\le\nu$ in a representation of highest weight $\nu$,
the Wigner functions \equ(WignerDef)
extend to entire analytic functions, and we have the bound
$$
\bigl|V^m_n(\cos(\beta/2),\vartheta)\bigr|
\le e^{|{\rm Im}(\beta/2)|n} e^{|{\rm Im}(\vartheta)||m|}\,,\qquad
\beta,\vartheta\in\complex\,.
\equation(VmnBound)
$$

The generalized spherical harmonics are known to be directly
related to the Jacobi polynomials $P^{(a,b)}_l$.
In particular,
$V^m_{m+2l}(r,\vartheta)=z^m P^{(0,m)}_l\bigl(2|z|^2-1\bigr)$
for $m,l\ge 0$, where $z=re^{i\vartheta}$.
Using the Rodrigues formula for the Jacobi polynomials,
one finds the following Rodrigues formula for the Zernike functions [\rKoo].
Let $k,l\ge 0$. Then
$$
V^m_n(r,\vartheta)=\VV^m_n\bigl(re^{i\vartheta},re^{-i\vartheta}\bigr)\,,\qquad
n=k+l\,,\quad m=k-l\,,
\equation(VVmn)
$$
where
$$
\VV^m_n=\partial_z^l\partial_{\bar z}^k\PP_n\,,\qquad
\PP_n(z,\bar z)={1\over n!}(z\bar z-1\bigr)^n\,.
\equation(VRodrigues)
$$
Here $\partial_z$ and $\partial_{\bar z}$ denote the partial
derivatives with respect to the
(independent) variables $z$ and $\bar z$, respectively.
The identity \equ(VRodrigues) can be used e.g.~to give a simple proof
of the following lemma.
Denote by $\Delta$ the Dirichlet Laplacean on the disk $\disk$.

\claim Lemma(InvLapV) {\rm [\rJan]}
Let $n\ge m\ge 0$ with $n-m$ even.
Then
$$
\eqalign{
\Delta^{-1}V^m_n&=c_2 V^m_{n+2}+c_1 V^m_n+c_0 V^m_{n-2}
\hskip20pt{\rm~if~} n>m\,,\cr
\Delta^{-1}V^m_n&=c_2V^m_{n+2}-c_2V^m_n
\hskip62pt{\rm~if~} n=m\,,\cr}
\equation(InvLapV)
$$
where
$$
c_2={1\over 4(n+2)(n+1)}\,,\quad
c_1=-{1\over 2n(n+2)}\,,\quad
c_0={1\over 4n(n+1)}\,.
\equation(InvLapConst)
$$

\proof
A trivial computation shows that
$$
\partial_z\partial_{\bar z}\PP_k=k\PP_{k-1}+\PP_{k-2}\,,\qquad k\ge 2\,.
\equation(ddbarPk)
$$
We claim that there exist constants $c_2$, $c_1$, and $c_0$,
such that
$$
c_2\bigl(\partial_z\partial_{\bar z}\bigr)^2\PP_{n+2}
+c_1\bigl(\partial_z\partial_{\bar z}\bigr)\PP_n
+c_0\PP_{n-2}=\tquarter\PP_n\,.
\equation(PreInvLap)
$$
Indeed, each term on the left hand side of this equation
is a linear combination of $\PP_n$, $\PP_{n-1}$, and $\PP_{n-2}$.
So we have $3$ linear equations with $3$ unknowns.
A straightforward computation yields the solution \equ(InvLapConst).
Applying $4\partial_z^l\bar\partial_{\bar z}^k$
to both sides of the equation \equ(PreInvLap), we obtain
$$
4\partial_z\partial_{\bar z}\bigl[c_2\VV^m_{n+2}+c_1\VV^m_n+c_0\VV^m_{n-2}\bigr]=\VV^m_n\,.
\equation(InvLap)
$$
Notice that $c_0+c_1+c_2=0$.
Thus, the function $[\ldots]$ in the equation \equ(InvLap) vanishes for $z\bar z=1$.
Taking $z=x+iy$ and $\bar z=x-iy$,
we have $4\partial_z\partial_{\bar z}=\partial_x^2+\partial_y^2=\Delta$,
and the first identity in \equ(InvLapV) follows.
The second identity is verified similarly.
\qed

To conclude this section,
we note that the change of variables $r=\cos(\beta/2)$
used in the definition \equ(VRfromD) is far from ad hoc.
It defines the azimuthal projection
$(\beta,\vartheta)\mapsto(r,\vartheta)$ from the sphere to the disk.
This projection preserves area (up to a trivial factor).
Using the orthogonality properties of the generalized
spherical harmonics, one finds that the Zernike functions $V^m_n$
constitute a complete orthogonal set for $\rmL^2(\disk)$.
Alternatively, one can use the orthogonality properties
of the Jacobi polynomials [\rWue].

\section Real analytic functions on the disk

In this section we prove that $\BB_\rho$
is a Banach algebra under pointwise multiplication of functions.
In addition, we give a bound on the inverse Dirichlet
Laplacean, and we introduce some notation that will be needed later on.
Unless specified otherwise, the domain of a function $u=u(r,\vartheta)$
is assumed to be the cylinder $[0,1]\times\sphere^1$.
But we still regard $u$ as a function on the disk $\disk$.

For every integer $m$, define $N_m=\{|m|,|m|+2,|m|+4,\ldots\}$.
Given a real number $\rho\ge 1$,
denote by $\AA_\rho$ the real vector space of all functions $u$,
$$
u=\sum_{m\in\integer}u_m\,,\qquad
u_m=\sum_{n\in N_m}u_{m,n}V^m_n\,,\qquad u_{m,n}\in\complex\,,
\equation(AArhoFun)
$$
that have a finite norm
$$
\|u\|_\rho=\sum_{m\in\integer}\|u_m\|_\rho\,,\qquad
\|u_m\|_\rho=\sum_{n\in N_m}|u_{m,n}|_{\sss 1}\,\rho^n\,.
\equation(AArhoNorm)
$$
Here $|x+iy|_{\sss 1}=|x|+|y|$ for $x,y\in\real$.
When equipped with this norm, $\AA_\rho$ is a Banach space over $\real$.
A real-valued function $u\in\AA_\rho$ can be written in the form
$$
u=\sum_{n\in N_0}A_{0,n}V^0_n
+\half\sum_{m\ne 0}\sum_{n\in N_m}\bigl[A_{m,n}-iB_{m,n}\bigr]V^m_n\,,
\equation(RealAArhoFun)
$$
with real coefficients $A_{m,n}$ and $B_{m,n}$
satisfying $A_{-m,n}=A_{m,n}$ and $B_{-m,n}=-B_{m,n}$, respectively,
for all integers $m$ and all $n\in N_m$.
Let us relabel the coefficients by setting
$a_{m,l}=A_{m,n}$ and $b_{m,l}=B_{m,n}$, where $l=(n-|m|)/2$.
Then a short computation shows that $u$
agrees with the function \equ(uZernikeSeries),
and that the norm \equ(AArhoNorm) of $u$ is given by \equ(BBrhoNormCS).
In other words,
$\BB_\rho$ is the subspace of real-valued functions $u\in\AA_\rho$.

\claim Lemma(BanchAlg)
$\AA_\rho$ is a Banach algebra under pointwise multiplication.
If $\rho>1$ then the functions in $\AA_\rho$ extend analytically
to some complex open neighborhood of $\disk$.

\proof
Consider first fixed integers $m_1$, $m_2$, and define $m_3=m_1+m_2$.
Recall from \equ(VProdRule) that
$$
V_{n_1}^{m_1}V_{n_2}^{m_2}
=\sum_{n_3}C_{n_1,n_2,n_3}^{m_1,m_2,m_3}\,V_{n_3}^{m_3}\,,\qquad
n_1\in N_{m_1}\,,\quad n_2\in N_{m_2}\,,
\equation(VProdRuleAgain)
$$
where $C_{n_1,n_2,n_3}^{m_1,m_2,m_3}$ is the square of a Clebsch-Gordan
coefficient and thus nonnegative.
These coefficients vanish whenever $n_j\not\in N_{m_j}$ for some $j$.
They also vanish if $n_3>n_1+n_2$, as we will describe later.
In addition, we have $\sum_{n_3}C_{n_1,n_2,n_3}^{m_1,m_2,m_3}=1$.
This follows from unitarity, but it can be seen also from \equ(VProdRuleAgain)
by noting that $R^m_n(1)=1$ whenever $n\in N_m$.
As a result,
$$
\|V_{n_1}^{m_1}V_{n_2}^{m_2}\|_\rho
\le\sum_{n_3}C_{n_1,n_2,n_3}^{m_1,m_2,m_3}\rho^{n_3}\le\rho^{n_1+n_2}\,.
\equation(VProdBound)
$$

Let now $u$ and $v$ be two functions in $\AA_\rho$.
To simplify notation, we define $u_{m,n}=0$
and $v_{m,n}=0$ whenever $n\not\in N_m$.
By using the bound \equ(VProdBound), we immediately get
$$
\eqalign{
\|uv\|_\rho
&\le\sum_{m_1,n_1,m_2,n_2}
|u_{m_1,n_1}v_{m_2,n_2}|_{\sss 1}\|V_{n_1}^{m_1}V_{n_2}^{m_2}\|_\rho\cr
&\le\sum_{m_1,n_1,m_2,n_2}|u_{m_1,n_1}|_{\sss 1}|v_{m_2,n_2}|_{\sss 1}
\rho^{n_1+n_2}
=\|u\|_\rho\|v\|_\rho\,.\cr}
\equation(BanachAlg)
$$
This shows that $\AA_\rho$ is a Banach algebra, as claimed.

Consider now $\rho>1$.
{}From the bound \equ(VmnBound), is is clear that
a function $u\in\AA_\rho$ extends analytically
to a complex open neighborhood $A$ of $[0,1]\times\sphere^1$
in the variables $(r,\vartheta)$.
So the series \equ(AArhoFun) for $u$
converges uniformly on compact subsets of $A$.
Changing to Cartesian variables $(x,y)$,
this translates into uniform convergence on
compact subsets of 
some open neighborhood $\disk_\rho$ of $\disk$.
But the Zernike functions $V^m_n$ are polynomials in $(x,y)$,
as can be seen e.g.~from \equ(VVmn).
Thus, being a locally uniform limit of analytic functions,
$u$ is analytic on $\disk_\rho$.
\qed

We note that Banach algebras of disk polynomials
have been considered before in [\rKanjin].

\smallskip
For our computer-assisted error estimates,
we approximate a function by truncating its Zernike series.
Given $N\ge 0$, define the projection $\proj_\ssN:\AA_\rho\to\AA_\rho$ as follows.
For every $u\in\AA_\rho$,
the function $\proj_\ssN u$ is obtained from $u$
by truncating the Zernike series \equ(AArhoFun) of $u$
to terms with index $n<N$.

\claim Proposition(InvLapBound)
Consider the inverse Dirichlet Laplacean $\Delta^{-1}$
as a linear operator on $\AA_\rho$.
Then $\Delta^{-1}$ is compact, and for every $N\ge 1$,
the operator norm of $\Delta^{-1}(\Id-\proj_\ssN)$ is bounded by
$$
\bigl\|\Delta^{-1}(\Id-\proj_\ssN)\bigr\|\le
{\bigl(\rho+\rho^{-1}\bigr)^2\over 4N(N+2)}\,.
\equation(InvLapBound)
$$

The bound \equ(InvLapBound) is an immediate consequence of \clm(InvLapV).
It implies in particular that $\Delta^{-1}$
is a uniform limit of finite rank operators, and thus compact.

\section Main steps in the proof

In this section we describe how
Theorems \clmno(HillRing), \clmno(HenonOne), and \clmno(HenonTwo)
can be proved by verifying the assumptions of five technical lemmas.
The estimates that are used to verify these assumptions
will be discussed in Section 6.

As mentioned in the introduction,
we find solutions of the equation \equ(Main)
by solving the fixed point problem for the map $G$ given by \equ(MainFix).
Here we follow the approach used in [\rAKi].
We always assume that $f'(u)=u^3$.
Let $\rho$ be a real number larger than $1$, to be specified later.
Assuming that $w$ belongs to $\BB_\rho$,
$G$ defines a smooth compact map on $\BB_\rho$.
This follows from the fact that $\BB_\rho$ is a Banach algebra, and from \clm(InvLapBound).

In this paper, we are interested only in solutions
$u=u(r,\vartheta)$ that are even functions of $\vartheta$.
The even subspace of $\BB_\rho$ will be denoted by $\BB_\rho^\even$.
Given $r>0$ and $u\in\BB_\rho^\even$,
denote by $B_r(u)$ the close ball in $\BB_\rho^\even$
of radius $r$, centered at $u$.

\smallskip
Given a function $\uzero\in\BB_\rho^\even$
and a bounded linear operator $M$ on $\BB_\rho^\even$, define
$$
\NN(h)=G(\uzero+Ah)-\uzero+Mh\,,\qquad  A=\Id-M\,,
\equation(NNDef)
$$
for every $h\in\BB_\rho^\even$.
Clearly, if $h$ is a fixed point of $\NN$ then 
$\uzero+Ah$ is a fixed point of $G$.
Furthermore, if the operator $\Id-DG(\uzero)$ is invertible,
and if $A$ sufficiently close to its inverse,
then $\NN$ is a contraction near $\uzero$.

Our goal is to apply the contraction mapping theorem
to the map $\NN$, on some small ball $B_r(0)$.
Thus $\uzero$ is chosen to be an approximate fixed point of $G$.
For practical reasons, $\uzero=\proj_\ssN \uzero$ for some $N$.
For the same reasons,
we choose $M$ to satisfy $M=\proj_\ssN M\proj_\ssN$ for some $N$.
So $M$ is in essence a matrix.

To guarantee the existence of a true fixed point of $G$ near $\uzero$,
it suffices to prove the hypotheses of the following lemma.

\claim Lemma(NNfix)
Let $\rho>1$ and $w\in\BB_\rho^\even$.
Assume that there exists a function $\uzero\in\BB_\rho^\even$,
a finite-rank operator $M:\BB_\rho^\even\to\BB_\rho^\even$,
and a real number $\delta>0$, such that
the map map $\NN$ defined by \equ(NNDef) admits bounds
$$
\eps\ge\|\NN(0)\|_\rho\,,\qquad
K\ge\|D\NN(h)\|\,,\quad\forall h\in B_\delta(0)\,,
\equation(NNfix)
$$
with $\eps$ and $K$ satisfying $\eps+K\delta<\delta$.
Then the equation \equ(Main) has a solution $u_\ast\in\BB_\rho$
within a distance $\|A\|\delta$ of $\uzero$.
Furthermore, if $M$ has no eigenvalue $1$,
then this solution $u_\ast$ is locally unique.

\proof
By the contraction mapping principle,
the given bounds imply that $\NN$ has a unique fixed point $h_\ast$
in the ball $B_\delta(0)$.
In fact, $h_\ast$ lies in the interior of $B_\delta(0)$,
since the inequality $\eps+K\delta<\delta$ is strict.
Clearly $u_\ast=\uzero+Ah_\ast$ is a fixed point of $G$.
If $M$ has no eigenvalue $1$, then this fixed point is locally unique,
since the fixed point $h_\ast$ of $\NN$
is locally unique and $A=\Id-M$ is invertible.
The distance of $u_\ast$ from $\uzero$ is
$\|u_\ast-\uzero\|_\rho=\|Ah_\ast\|_\rho\le\|A\|\delta$.
\qed

With $u_\ast$ as above,
consider the possibility that $|u_\ast|$ is invariant
under some nontrivial rotation.
Then the function $u=u_\ast^2$ is invariant under a rotation by $2\pi/k$
for some integer $k\ne 1$.
So the component $u_m$ in the representation \equ(AArhoFun) vanishes,
unless $m$ is a multiple of $k$.
This proves the following.

\claim Lemma(MinSymm)
If some coefficient $u_{1,l}$
in the Zernike expansion \equ(uZernikeSeries) for $u=u_\ast^2$ is nonzero,
then $|u_\ast|$ is not invariant under any nontrivial rotation.

This fact is used to verify that the solutions
described in Theorems \clmno(HillRing) and \clmno(HenonOne)
have no rotation symmetries.
The property $u_\alpha>0$ mentioned in \clm(HenonOne) follows
from the well-known fact that index-$1$ solutions
do not vanish anywhere on $\disk$, if $w>0$.
In order to prove the symmetry properties
of the solutions described in \clm(HenonTwo), we use the following.
We assume that $w$ is radial.

\claim Lemma(HasSymm)
Under the assumptions of \clm(NNfix),
if $\uzero$ is invariant under $S_n$,
and if $M$ commutes with $S_n$, then the solution $u_\ast$
described in (the proof of) \clm(NNfix) is invariant under $S_n$.

This claim follows from the fact that,
under the given assumptions, if $h\in\BB_\rho^\even$
is invariant under $S_n$, then so is $\NN(h)$.
Thus, the limit $\displaystyle h_\ast=\lim_{k\to\infty}\NN^k(0)$
and the function $u_\ast=\uzero+h_\ast-Mh_\ast$ are invariant under $S_n$.

\smallskip
Next, we consider the problem of determining
the Morse index of a solution $u$ of the equation \equ(Main).
As in [\rAKi] we use the identity
$$
D^2J(u)(v\times v)
=\int_\disk\,\Bigl[
|\nabla v|^2-3wu^2v^2\Bigr]\,dxdy
=\Langle v,[\id-DG(u)]v\Rangle_{\rmH^1}\,,
\equation(DDJ)
$$
which relates the second derivative of the functional $J$ defined in \equ(Ju)
to the first derivative of the map $G$ defined in \equ(MainFix).
Here, it is assumed that $v$ belongs to $\rmH^1_0=\rmH^1_0(\disk)$.
Notice that, if $W$ is a bounded linear operator on $\rmL^2=\rmL^2(\disk)$,
then $(-\Delta)^{-1}W$ is a bounded linear operator
from $\rmL^2$ to $\rmH^1_0$, and
$$
\Langle(-\Delta)^{-1}Wv,h\Rangle_{\rmH^1}
=\langle Wv,h\rangle_{\rmL^2}\,,\qquad
v,h\in\rmH^1_0\,.
\equation(HiLiiProds)
$$
Clearly, if $W$ is self-adjoint on $\rmL^2$,
then the restriction of $(-\Delta)^{-1}W$ to $\rmH^1_0$ is self-adjoint.

These observations explain much of the following.

\claim Proposition(Morse)
Assume that $wu^2$ is of class $\rmC^1$ and nonnegative.
Then $DG(u)$ is a compact positive self-adjoint operator on $\rmH^1_0$.
If $wu^2$ is positive almost everywhere on $\disk$,
then all eigenvalues of $DG(u)$ are positive.
Assume now that $wu^2$ belong to $\BB_\rho$.
Then every eigenvector of $DG(u)$ with nonzero eigenvalue
belongs to $\BB_\rho$.
If in addition $u$ solves the equation \equ(Main),
then the Morse index of $u$ equals
the number of eigenvalues of $DG(u)$ that exceed $1$.

This proposition was proved in [\rAKi] for a square domain.
The same arguments apply in the case of a disk,
using that $\Delta^{-1}$ defines a compact linear operator on $\BB_\rho$ by \clm(InvLapBound),
and that $\BB_\rho$ is dense in $\rmH^1$.
The density of $\BB_\rho$ in $\rmH^1$ follows e.g.~from the fact
that the Zernike functions $V^m_n$
constitute a complete orthogonal set for $\rmL^2$.

\smallskip
Assume now that $u\in\BB_\rho$ is a nontrivial fixed point of $G$, with $w\in\BB_\rho$.
By \clm(Morse), the Morse index of $u$
agrees with the number of eigenvalues of $DG(u)$ that are larger than $1$.
And it suffices to consider $DG(u)$ as a linear operator on $\BB_\rho$.

Notice that $\BB_\rho^\even$ is an invariant subspace of $DG(u)$.
Another invariant subspace of $DG(u)$ is the space $\BB_\rho^\odd$
of all functions $u=u(r,\vartheta)$ in $\BB_\rho$ that are odd functions of $\vartheta$.
We refer to $\BB_\rho^\odd$ as the odd subspace of $\BB_\rho$.
Clearly, every eigenvalue of $DG(u)$ has an eigenfunction
that belongs to one of these two subspaces.
Two eigenvalues are known explicitly:
$u\in\BB_\rho^\even$ is an eigenvector of $DG(u)$ with eigenvalue $3$,
and $\partial_\vartheta u\in\BB_\rho^\odd$ is an eigenvector of $DG(u)$ with eigenvalue $1$.
This follows from the fact that nonlinearity $f'(u)$ in the equation \equ(MainFix)
is cubic, and that $(r,\vartheta)\mapsto u(r,\vartheta+t)$
is a fixed point of $G$ for every real number $t$, respectively.

Consider now the restriction of $DG(u)$
to one of the subspaces $\BB_\rho^\parity$ of fixed parity $\paritysym\in\{\evensym,\oddsym\}$.
Our goal is to determine the number of eigenvalues of $DG(u)$
that are larger than $1$.
In order simplify our description,
let $\lambda_1\ge\lambda_2\ge\ldots\ge 0$ be the eigenvalues of $DG(u)$,
listed with their multiplicities,
and let $u_1,u_2,\ldots$ be the corresponding eigenvectors.
We may assume that these eigenvectors are mutually orthogonal.
If $\lambda_n<\theta<\lambda_{n+1}$, then the operator
$$
DG(u)-\KK\,,\qquad
\KK=\sum_{j=1}^n\lambda_j
{\langle h,u_j\rangle_{\sss\rmH^1}
\over\langle u_j,u_j\rangle_{\sss\rmH^1}}u_j\,,
\equation(DGuMinusKK)
$$
has a spectral radius less than $\theta$.
We are interested in obtaining a similar conclusion
by using only approximate eigenvalues and eigenvectors.

The is possible by using the following fact.

\claim Lemma(atmostn)
Let $A$ and $K$ be bounded linear operators on a Hilbert space.
Assume that $A$ is normal, that $K$ is of finite rank $n$,
and that $\|A-K\|<\theta$.
Then $A$ has at most $n$
eigenvalues $\lambda_j$ (counting multiplicities) satisfying $|\lambda_j|\ge\theta$.

\proof
As a rank $n$ operator, $K$ admits a representation
$
Kh=\sum_{i=1}^n a_i\langle h,w_i\rangle v_i\,.
$
Assume for contradiction that $A$ admits
an orthonormal set $\{u_1,u_2,\ldots,u_{n+1}\}$
of $n+1$ eigenvectors with eigenvalues $\lambda_1,\lambda_2,\ldots\lambda_{n+1}$
satisfying $|\lambda_j|\ge\theta$.
Then some nontrivial linear combination
$h=\sum_{j=1}^{n+1}c_ju_j$ is orthogonal to each of the $n$ vectors $w_i$.
It satisfies $Kh=0$, and thus
$$
\bigl\|(A-K)h\bigr\|^2
=\|Ah\|^2=\sum_{j=1}^{n+1}|\lambda_j c_j|^2
\ge\theta\sum_{j=1}^{n+1}|c_j|^2=\theta\|h\|^2\,.
\equation(preqed)
$$
This contradicts the assumption that $\|A-K\|<\theta$.
\qed

In our application, $\theta<0.981$,
and $A$ is the restriction of $DG(u)$ to the subspace $\BB_\rho^\parity$.
To be more specific:
for the solutions $u=u_w,u_2,u_4,u_6,u_{2,2},u_{2,4}$
described in Theorems \clmno(HillRing), \clmno(HenonOne), and \clmno(HenonTwo),
we verify the assumptions of this proposition
with $n=2,1,1,1,2,3$, respectively, on the even subspace,
and with $n=1,1,1,1,1,2$, respectively, on the odd subspace.

\clm(atmostn) shows that $n$ is an upper bound
on the number of eigenvalues of $A$ in the interval $[\theta,\infty)$.
We need a lower bound as well; but only if $n>1$,
since $A$ has an eigenvalue $3$ or $1$, as described earlier.
Such a bound can obtained by using the following lemma,
for some real number $a>1$.

\claim Lemma(atleastm)
Let $A$ be a compact self-adjoint linear operator on a Hilbert space $\HH$.
Let $\{v_1,v_2,\ldots,v_m\}$
be an orthonormal set in $\HH$, and assume that
$$
A_{j,j}-\sum_{i\ne j}|A_{i,j}|>a\,,\qquad
A_{i,j}=\langle v_i,A v_j\rangle\,,\quad 1\le i,j\le m\,.
\equation(atleastm)
$$
Then $A$ has at least $m$
eigenvalues (counting multiplicities) in the interval $[a,\infty)$.

This lemma is an immediate consequence
of the Gershgorin circle theorem and Cour\-ant's min-max principle.
We verify the hypotheses with $m=2,1,1,1,2,3$, respectively,
for the solutions $u=u_w,u_2,u_4,u_6,u_{2,2},u_{2,4}$ on the even subspace,
and with $m=1$ for the solution $u=u_{2,4}$ on the odd subspace.
The value $a=3$ works in all cases.
More accurate eigenvalue bounds can be found in [\rFiles].

We remark that this method for determining the Morse index
is significantly simpler,
and more efficient, than the method used in [\rAKi].

\section Clebsch-Gordan coefficients

In this section we give a brief description of
the identities and algorithms that we use to compute
and index Clebsch-Gordan coefficients.
For details we refer to the Ada packages {\tt Regges} and {\tt CG} in [\rFiles].

The Clebsch-Gordan coefficients 
$\langle\nu_1\mu_1\nu_2\mu_2|\nu_3\mu_3\rangle$
that appear in the product expansion \equ(VProdRule) vanish
unless the angular momenta $\nu_j$ and $\mu_j$ satisfy certain constraints.
These constraints, as well as symmetries,
are most conveniently described in terms of the Regge symbol
$$
\boxit{2pt}{
 $\matrix{
 \nu_2\!+\!\nu_3\!-\!\nu_1 & \nu_3\!+\!\nu_1\!-\!\nu_2 & \nu_1\!+\!\nu_2\!-\!\nu_3\cr
 \nu_1-\mu_1       & \nu_2-\mu_2       & \nu_3-\mu_3\cr
 \nu_1+\mu_1       & \nu_2+\mu_2       & \nu_3+\mu_3\cr}$}
={(-1)^{\nu_1-\nu_2+\mu_3}\over\sqrt{2\nu_3+1}}
\Langle\nu_1\mu_1\nu_2\mu_2|\nu_3(-\mu_3)\Rangle\,.
\equation(ReggeSym)
$$
A Regge symbol $\boxit{1.5pt}{$R$}$
vanishes unless its entries $R_{ij}$ are all nonnegative integers,
and unless the row sums $R_{i1}+R_{i2}+R_{i3}$
and column sums $R_{1j}+R_{2j}+R_{3j}$ all have the same value.
Furthermore, its absolute value is invariant under interchanges of rows,
interchanges of columns, and transposition [\rRegge].
Under odd row or column permutations of the matrix $R$,
the symbol acquires a factor $(-1)^J$,
where $J=\nu_1+\nu_2+\nu_3$.

A commonly used formula [\rRBMW] for Wigner's $3j$ symbol,
expressed here in terms of the Regge symbol, is
$$
\boxit{1.5pt}{$R$}=(-1)^{R_{12}-R_{33}}
\sqrt{{R_{11}!R_{12}!\cdots R_{2,3}!R_{3,3}!\over(J+1)!}}\sum_k{(-1)^k\over Q_k(R)}\,,
\equation(RBMW)
$$
where
$$
Q_k(R)=k!(R_{12}-R_{21}+k)!
(R_{11}-R_{32}+k)!(R_{31}-k)!(R_{21}-k)!(R_{32}-k)!\,.
\equation(QkR)
$$
The sum in \equ(RBMW) runs over all integers $k$
such that the arguments of all factorials in \equ(QkR) are nonnegative.

As can be seen from \equ(RBMW),
the square of a Regge symbol is a rational number $P/Q$.
In our programs, we compute $P$ and $Q$ exactly,
following roughly a procedure described in [\rJoFo].
The summands $(-1)^k/Q_k$ in the equation \equ(RBMW)
are multiplied first by their least common multiple,
so the sum becomes a sum of integers.
Factorials and their products
are computed in terms of their prime factorization.

This computation is too costly to be repeated
whenever we need the value of a Regge symbol.
Thus, we compute the necessary values
beforehand and store them in a linear array.
Due to the above-mentioned symmetries,
it suffices to index Regge matrices of the form
$$
\boxit{1.5pt}{$\RR$}=\boxit{2pt}{
 $\matrix{
     S &     L & X+B-T\cr
     X &     B & S+L-T\cr
 L+B-T & S+X-T &     T\cr}$
}\,,
\equation(ReggePerm)
$$
with $L\ge X\ge T\ge B\ge S\ge 0$.
As was shown in [\rRaYu],
any Regge matrix $R$ with nonzero symbol
can be transformed to such a normal form $\RR$
via row and column permutations, and possibly a transposition.
This eliminates half of the $72$ symmetries.

Here, we eliminate much of the other half as well by requiring that $T\le(L+S)/2$.
This can be achieved by exchanging the last two rows of $\RR$, if necessary.

We index such Regge matrices by enumerating the set
$$
\SS=\bigl\{(l,s,t,x,b):\;\;
l\ge x\ge t\ge b\ge s\ge 0\,,\;\;
t\le (l+s)/2\,\bigr\}\,,
\equation(OrderedIndices)
$$
using the following lexicographical order:
recursively, define $(u,v,w,\ldots)<(U,V,W,\ldots)$
to mean that either $u<U$, or else $u=U$ and $(v,w,\ldots)<(V,W,\ldots)$.
Assuming that the quintuple $(L,S,T,X,B)$ belongs to $\SS$,
the index of the Regge matrix \equ(ReggePerm)
is defined to be the number of quintuples in $\SS$
that are less than or equal to $(L,S,T,X,B)$.
This index can be expressed as
$$
\II(\RR)=\sum_{l\le L}\sum_{s\le S}\sum_{t<T}\sum_{x,b}\Chi_\SS(l,s,t,x,b)
+C_{S,T,X,B}\,,
\equation(IICCp)
$$
where $C_{S,T,X,B}$ is the sum over all $x\le X$ and $b\le B$
of the numbers $\Chi_\SS(L,S,T,x,b)$.
Here $\Chi_\SS$ denotes the indicator function of the set $\SS$.
Notice that the five-fold sum in \equ(IICCp)
defines a function of the three variables $(L,S,T)$.
In our programs, the values of this function are determined simply by having the computer
carry out the five-fold sum.
The values are then stored in a three-dimensional array.
As for the values $C_{S,T,X,B}$,
a straightforward computation shows that
$$
C_{S,T,X,B}=(X-T)(T-S+1)+(B-S+1)\,.
\equation(CpVal)
$$

\section Computer estimates

In order to complete our proof of
Theorems \clmno(HillRing), \clmno(HenonOne), and \clmno(HenonTwo),
we need to verify the assumptions of the lemmas in Section 4.
This is done with the aid of a computer.
For each of the six models considered, we have chosen $\rho=65/64$.

To fix ideas, consider \clm(NNfix) for some given choice of $w\in\BB_\rho$.
To verify the assumptions of this lemma,
we first determine an approximate fixed point $\uzero$ of $G$
and an approximation $M$ for the operator $\Id-[\Id-DG(\uzero)]^{-1}$.
These numerical data are included with the source code of our programs in [\rFiles].
The remaining steps are rigorous:
First, we compute an upper bound $\eps$ on the norm of $\NN(0)$.
Using this bound, we define an increasing function $d$ on the interval $[0,3/4]$
such that $d(K)>\eps/(1-K)$ on this interval.
Now we compute an upper bound $K$ on the operator norm of $D\NN(h)$
that holds for all $h$ of norm $d(3/4)$ or less.
After verifying that $K\le 3/4$, we set $\delta=d(K)$.
This guarantees that $\eps+K\delta<\delta$.
We also verify that $M$ has no eigenvalue $1$.

The rigorous part is still numerical,
but instead of truncating series and ignoring rounding errors,
it produces guaranteed enclosures
at every step along the computation.
Our choice of enclosures and associated data types will be described below.

The above-mentioned steps are analogous to
those used in the proof of Theorem 4.1 in [\rAKi].
The main difference is that [\rAKi]
uses functions on the square and data of type {\tt Fourier2},
while here we use functions on the disk and data of type {\tt Zernike}.
To avoid undue repetition,
we will focus here on those aspects of the proof
where the differences are relevant.

We will also describe our computation of the Morse index,
which amounts to verifying the assumptions of the Lemmas \clmno(atmostn) and \clmno(atleastm).
But any description given here is necessarily incomplete.
For precise definitions and other details,
the ultimate reference is the source code of our programs [\rFiles].
This code is written in the programming language Ada [\rAda].

One of the basic data type in our programs
is the type {\tt Ball} that we use to define enclosures for real numbers.
A data item of type {\tt Ball} is a pair {\tt B=(B.C,B.R)},
where {\tt B.C} is a representable number (type {\tt Rep}),
and where  {\tt B.R} a nonnegative representable number (type {\tt Radius}).
The corresponding subset of $\real$ is the interval
${\tt B}^\flat=\{b\in\real: |b-{\tt B.C}|\le{\tt B.R}\}$.
Using controlled rounding, it is trivial to implement e.g.~a
``{\tt function Sum(A,B: Ball) return Ball}''
with the property that ${\tt Sum(A,B)}^\flat$ contains $a+b$ whenever
$a\in{\tt A}^\flat$ and $b\in{\tt B}^\flat$.
Similarly for other elementary operations involving real numbers.

Next, we describe our enclosures for functions in $\BB_\rho$
that belong to the even subspace $\BB_\rho^\even$
or to the odd subspace $\BB_\rho^\odd$.
The enclosures depend on the choice of a positive integer {\tt Size}
which we denote here by $S$. Define $D=\lfloor S/2\rfloor$.
We start by considering functions $f:[0,1]\to\real$
with the property that
$$
(J_\rho^{\pm m}f)(r,\vartheta)\defeq f(r)e^{\pm im\vartheta}
\equation(Jmrhof)
$$
defines a function in $\BB_\rho$.
In this step, $\rho\ge 1$ and $m\ge 0$ are considered fixed,
with $m\le D$.
Our enclosures for such functions $f$ are associated
with a data type {\tt Radial}.
A data item of type {\tt Radial} is (in essence) a triple {\tt F=(F.M,F.C,F.E)},
where {\tt F.C} is an {\tt array(0 .. D) of Ball},
{\tt F.E} is an {\tt array(0 .. D+1) of Radius},
and ${\tt F.M}=m$.
The corresponding set ${\tt F}^\flat$
is the set of all function $f:[0,1]\to\real$ that admit a representation
$$
f(r)=\sum_{j=0}^{D_m}C_j R^m_{m+2j}(r)
+\sum_{j=0}^{D_m+1} E_j(r)\,,\qquad D_m=\lfloor(D-m)/2\rfloor\,,
\equation(fRadial)
$$
with $C_j\in{\tt F.C(j)}^\flat$ for $j\le D_m$,
and $\|J_\rho^m E_j\|_\rho\le{\tt F.E(j)}$ for $j\le D_m+1$.
In addition, we require that the Zernike series
for the functions $J_\rho^m E_j$ include only modes $V^m_{n'}$
with $n'\ge n$, where $n=m+2j$.
Notice that the coefficient array {\tt F.C}
specifies a set of polynomials of degree $\le m+2D_m\le D$.
The numbers in {\tt F.E} represent error bounds.

An item of type {\tt Zernike} is a quadruple {\tt U=(U.R,U.P,U.C,U.E)},
where ${\tt U.R}\ge 1$ is of type {\tt Radius},
{\tt U.C} is an {\tt array(0 .. S) of Radial}
with {\tt U.C(m).M=m} fixed for each {\tt m},
{\tt U.E} is an {\tt array(0 .. 2*S) of Radius},
and {\tt U.P} is either $0$ or $1$.
If ${\tt U.P}=1$ then {\tt U} defines a subset ${\tt U}^\flat$
of $\BB_\rho^\odd$ with $\rho={\tt U.R}$.
Consider now ${\tt U.P}=0$. In this case,
{\tt U} defines a subset ${\tt U}^\flat$ of $\BB_\rho^\even$ with $\rho={\tt U.R}$.
This set consists of all functions
$$
u=\sum_{m=0}^S \half\bigl(J_\rho^m+J_\rho^{-m}\bigr)f_m
+\sum_{m=0}^{2S} E_m\,,\qquad E_m\in\BB_\rho^\even\,,
\equation(uZernike)
$$
with $f_m\in{\tt U.C(m)}^\flat$ for $0\le m\le S$,
and $\|E_m\|_\rho\le{\tt U.E(m)}$ for $0\le m\le 2S$.
In addition, we require that the Zernike series for $E_m$ include
only modes $V^{m'}_n$ with $m'\ge m$.
Our enclosures for $\BB_\rho^\odd$ are defined analogously.

{\tt Zernike}-type sets ${\tt U}^\flat$
play the same role for functions in $\BB_\rho^\even\cup\BB_\rho^\odd$
as {\tt Ball}-type sets ${\tt B}^\flat$ play for real numbers.
It is trivial to implement e.g.~a ``{\tt function Sum(U,V: Zernike) return Zernike}''
with the property that ${\tt Sum(U,V)}^\flat$ contains $u+v$ whenever
$u\in{\tt U}^\flat$ and $v\in{\tt V}^\flat$, provided that {\tt U.P=V.P}.
Implementing a bound on the product of two such functions
is a bit more involved.
Here we use the Banach algebra property of $\BB_\rho$
and enclosures for the Clebsch-Gordan coefficients.
For details we refer to the Ada package {\tt Zernikes} in [\rFiles].
This package also implements a bound {\tt InvNegLap}
on the operator $(-\Delta)^{-1}$,
using estimates of the type \equ(InvLapBound)
for the error terms in \equ(fRadial) and \equ(uZernike).

More problem-specific operations are defined in {\tt Zernikes.GFix},
including bounds {\tt GMap} and {\tt DGMap}
on the map $G$ and its derivative $DG(u)$, respectively.
Our proof of \clm(NNfix) is organized by the procedure {\tt ContrFix},
using {\tt DContrNorm} to obtain a bound on the operator norm of $D\NN(h)$.
The steps are as described at the beginning of this section.
This applies to each of the
solutions $u_w$, $u_\alpha$, and $u_{2,n}$, described in
Theorems \clmno(HillRing), \clmno(HenonOne), and \clmno(HenonTwo), respectively.
For the solutions $u_w$ and $u_\alpha$
we also verify the assumptions of \clm(MinSymm),
and for $u_{2,n}$ we verify the assumptions of \clm(HasSymm).
The details can be found in [\rFiles].

What remains to be discussed is the computation of the Morse index.
Using Lemmas \clmno(atmostn) and \clmno(atleastm),
with $A$ being the restriction of $DG(u)$ to $\BB_\rho^\sigma$,
this task is relatively straightforward.
The computations for $\sigma=0$ and for $\sigma=1$ are carried out separately.
And this is done for each of the six models being considered.

Among the data included in [\rFiles] are approximate eigenvectors of $A$.
They define a self-adjoint approximation $K$
for the operator $\KK$ described in \equ(DGuMinusKK).
A bound on the the map $\LL=A-K$ is implemented by the procedure {\tt LLMap}.
In order to estimate the spectral radius of $\LL$, as required by \clm(atmostn),
we first construct an enclosure {\tt L}
for the operator $\LL:\BB_\rho^\parity\to\BB_\rho^\parity$,
iterate ${\tt L}\mapsto{\tt L}^2$ several times,
and then estimate the operator norm of the result.
The inequalities \equ(atleastm) needed for \clm(atleastm)
are verified in the procedure {\tt KBound}.
This procedure first orthonormalizes the approximate
eigenvectors that were used to define the operator $K$.

In order to construct operator enclosures,
we use some data types and procedures from {\tt Zernikes}
that we have not yet described.
Notice that a {\tt Zernike} {\tt U}
can be viewed as a collection of ``coefficient modes'' {\tt U.C(m).C(j)}
and ``error modes'' {\tt U.C(m).E(j)} or {\tt U.E(m)}.
Coefficient modes represent one-dimensional
subspaces of $\BB_\rho^\parity$, while error modes
represent infinite-dimensional subspaces.
To specify individual modes we use a data type {\tt ZMode}.
We are interested mostly in finite collections of modes
whose subspaces $Z_i$ define a partition of $\BB_\rho^\parity$,
in the sense that $\bigoplus_i Z_i=\BB_\rho^\parity$,
and that $Z_i\cap Z_j=\{0\}$ for $i\ne j$.
Such a ``partition'' is specified by our data type {\tt ZModes}.
Our linear operator $\LL:\BB_\rho^\parity\to\BB_\rho^\parity$
now defines a ``matrix'' of operators $\LL_{i,j}:Z_j\to Z_i$.
By an enclosure for $\LL$ we mean a corresponding matrix of bounds,
with each element being a {\tt Ball}.
To be more precise, we restrict to {\tt ZModes} that
allow a {\tt Zernike} to be distributed efficiently
over the individual modes, using the procedure {\tt Extract}.
Then a bound ${\tt L}_{i,j}$ on $\LL_{i,j}$ is obtained
in essence by applying {\tt LLMap} to the $j$-th {\tt ZMode}
and extracting the $i$-th {\tt ZMode} from the result.

All major steps that are used to verify the assumptions
of the five lemmas in Section 4 are implemented in the procedures
described above.
They are combined in the proper order,
and invoked with the appropriate parameters,
by the main program {\tt Run\_All}.
Instructions on how to compile and run this program
are in a file {\tt README} that is included
with the source code of our programs in [\rFiles].
The programs {\tt Find\_Fix} and {\tt Find\_Eigen}
that were used to compute our numerical data are included as well.

The parameter {\tt Size} that determines the size of our {\tt Zernike}-type data
ranges from $70$ to $140$, depending on the computation.
For the set of representable numbers ({\tt Rep})
we choose standard extended floating-point numbers [\rIEEE]
that support controlled rounding,
and for bounds on non-elementary {\tt Rep}-operations
we use the open source MPFR library [\rMPFR].
Our programs were run successfully on a standard
desktop machine, using a public version of the gcc/gnat compiler [\rGnat].

\references

{\ninepoint


\item{[\rCNZ]} G.~Chen, W.-M. Ni, J.~Zhou,
{\sl Algorithms and visualization for solutions of nonlinear elliptic equations},
Int. J. Bifurc. Chaos {\bf 10}, 1565--1612 (2000).

\item{[\rSSW]} D.~Smets, J.~Su, M.~Willem,
{\sl Non radial ground states for the H\'enon equation},
Comm. Contemp. Math. {\bf 4} 467--480 (2002).

\item{[\rPW]} F.~Pacella, T.~Weth,
{\sl Symmetry of solutions to semilinear elliptic equations via Morse index},
Proc. Amer. Math. Soc. {\bf 135}, 1753--1762 (2007).

\item{[\rBaCa]} M.~Badiale, G.~Cappa,
{\sl Non radial solutions for non homogeneous H\'enon equation},
Nonlinear Analysis {\bf 109}, 45--55 (2014).


\item{[\rKoo]} T.~Koornwinder,
{\sl Two-variable analogues of the classical orthogonal polynomials},
in: Theory and application of special functions, R.A.~Askey (ed.),
Academic Press, 1975, pp.~435-495.

\item{[\rKin]} E.C.~Kintner,
{\sl On the mathematical properties of the Zernike Polynomials}.
Opt. Acta. {\bf 23} 679--680 (1976).

\item{[\rTango]} W.J.~Tango,
{\sl The Circle Polynomials of Zernike and Their Application in Optics},
Appl.~Phys. {\bf 13} 327--332 (1977).

\item{[\rKanjin]} Y.~Kanjin,
{\sl Banach algebra related to disk polynomials},
T\^ohoku Math. Journ. {\bf 37} 395--404 (1985).

\item{[\rWue]} A.~W\"unsche,
{\sl Generalized Zernike or disc polynomials},
J. Comput. Appl. Math. {\bf 174}, 135--163 (2005).

\item{[\rJan]} A.J.E.M.~Jansen,
{\sl Zernike expansion of derivatives and Laplacians
of the Zernike circle polynomials},
J. Opt. Soc. Am. A {\bf 31}, 1604--1613 (2014).


\item{[\rAKi]} G.~Arioli, H.~Koch,
{\sl Non-symmetric low-index solutions for a symmetric
boundary value problem},
J. Diff. Equations {\bf 252}, 448--458 (2012).

\item{[\rAKii]} G.~Arioli, H.~Koch,
{\sl Some symmetric boundary value problems and non-symmetric solutions},
J. Diff. Equations {\bf 259}, 796--816 (2015).

\item{[\rCZ]} J.~Cyranka, P.~Zgliczy\'nski,
{\sl Existence of globally attracting solutions for one-dimensional viscous Burgers equation
with nonautonomous forcing - a computer assisted proof},
SIAM J. Appl. Dyn. Syst. {\bf 14}, 787--821 (2015).

\item{[\rWaNa]}Y.~Watanabe, M.T.~Nakao,
{\sl A numerical verification method for nonlinear functional equations
based on infinite-dimensional Newton-like iteration},
Appl. Math. Comput. {\bf 276} 239--251 (2016).

\item{[\rCLJii]} R.~Castelli, J.-P.~Lessard, J.D.~Mireles-James,
{\sl Parameterization of invariant manifolds for periodic orbits (II):
A posteriori analysis and computer assisted error bounds},
J. Dyn. Diff. Equat. (2017).
{\tt https://doi.org/10.1007/s10884-017-9609-z}

\item{[\rFH]} J.L.~Figueras, A.~Haro,
{\sl Rigorous computer assisted application of KAM theory: a modern approach},
A. Found. Comput. Math. {\bf 17}, 1123--1193 (2017).

\item{[\rFL]} J.-L.~Figueras and R.~de la Llave,
{\sl Numerical computations and computer assisted proofs of periodic
orbits of the Kuramoto-Sivashinsky equation},
SIAM J. Appl. Dyn. Syst. {\bf 16}, 834--852 (2017).

\item{[\rBVCDLVWY]} I.~Bal\'azs, J.B.~van den Berg, J.~Courtois, J.~Dud\'as,
J.-P.~Lessard, A.~V\"or\"os-Kiss, J.F.~Willi\-ams, X.Y.~Yin,
{Computer-assisted proofs for radially symmetric solutions of PDEs},
preprint (2017)

\item{[\rAKiii]} G.~Arioli, H.~Koch,
{\sl Spectral stability for the wave equation with periodic forcing},
Preprint {\tt mp\_arc 17-23}.

\item{[\rMore]}
For earlier work, see references in [\rAKi\range\rAKiii].


\item{[\rRegge]} T.~Regge,
{\sl Symmetry properties of Clebsch-Gordan’s coefficients},
Nuovo Cimento {\bf 10},544--545 (1958).

\item{[\rRBMW]}M.~Rotenberg, R.~Bivins, N.~Metropolis, and J.K.~Wooten,
{\sl The 3j and 6j symbols}, Cambridge, MA: MIT Press, 1959.

\item{[\rRaYu]} J.~Rasch and A.C.H.~Yu,
{\sl Efficient storage scheme for precalculated Wigner
$3J$, $6J$ and Gaunt coefficients},
Siam J. Sci. Comput {\bf 25}, 1416--1428 (2003).

\item{[\rJoFo]} H.T.~Johansson and C.~Forss\'en,
{\sl Fast and accurate evaluation of Wigner $3j$, $6j$, and $9j$ symbols
using prime factorisation and multi-word integer arithmetic},
SIAM J. Sci. Comput. {\bf 38}, A376–A384 (2016) 


\item{[\rFiles]} G.~Arioli, H.~Koch,
The computer programs and data files are available at\hfill\break
{\tt http://www.ma.utexas.edu/users/koch/papers/zerni/}

\item{[\rAda]} Ada Reference Manual, ISO/IEC 8652:2012(E),
available e.g. at\hfil\break
{\tt http://www.ada-auth.org/arm.html}

\item{[\rGnat]} 
A free-software compiler for the Ada programming language,
which is part of the GNU Compiler Collection; see
{\tt http://gnu.org/software/gnat/}

\item{[\rIEEE]} The Institute of Electrical and Electronics Engineers, Inc.,
{\sl IEEE Standard for Binary Float\-ing--Point Arithmetic},
ANSI/IEEE Std 754--2008.

\item{[\rMPFR]} The MPFR library for multiple-precision floating-point computations
with correct rounding; see
{\tt http://www.mpfr.org/}.

}
\bye